\documentclass[preprint,10pt]{elsarticle}

\usepackage[T1]{fontenc}
\usepackage[utf8]{inputenc}
\usepackage[onehalfspacing]{setspace}
\usepackage[english]{babel}
\usepackage[left=3cm,right=2cm,top=3cm,bottom=3cm]{geometry}

\usepackage{tikz}
\usetikzlibrary{decorations.text, positioning, shapes.geometric, arrows}
\usepackage{xcolor}
\usepackage[english,vlined,algoruled, linesnumbered]{algorithm2e}
\usepackage{paralist}
\usepackage{ulem}
\usepackage{amsmath}
\usepackage{amssymb} 
\usepackage{tabularx, multirow}
\usepackage[flushleft]{threeparttable}
\usepackage{booktabs, caption}
\usepackage{marginnote}
\usepackage{longtable}
\usepackage{subcaption}
\usepackage{rotating}

\usepackage{makeidx}
\usepackage[nonumberlist,toc,acronym,nopostdot]{glossaries}

\usepackage{natbib}
 \bibpunct[, ]{(}{)}{,}{a}{}{,}%
\usepackage{graphicx}

\newcommand{\comment}[1]{}

\usepackage{graphicx} 
\usepackage{booktabs} 
\usepackage{xparse}   
\NewDocumentCommand{\rot}{O{45} O{1em} m}{\makebox[#2][l]{\rotatebox{#1}{#3}}}%
\usepackage{rotating}
\usepackage{floatrow}
\usepackage{pgfplots}
\DeclareUnicodeCharacter{2212}{−}
\usepgfplotslibrary{groupplots,dateplot}
\usetikzlibrary{patterns,shapes.arrows}
\pgfplotsset{compat=newest}

\newacronym[description=Adaptive Large Neighborhood Search]{ALNS}{ALNS}{adaptive large leighborhood search}

\newacronym[description=callbacks]{CB}{CB}{callback procedure}
\newacronym[description=constructive LB]{cLB}{cLB}{constructive lower bound}
\newacronym[description=construction heuristic]{CH}{CH}{construction heuristic}
\newacronym[description=destructive LB]{dLB}{dLB}{destructive lower bound}
\newacronym[description=column generation]{CG}{CG}{column generation}
\newacronym[description=constraint Programming]{CP}{CP}{constraint programming}
\newacronym[description=Composite Neighborhood Descent]{CND}{CND}{composite neighborhood descent}
\newacronym[description=crew scheduling]{CS}{CS}{crew scheduling}

\newacronym[description=Destructive Bound Improvement]{DBI}{DBI}{destructive bound improvement phase}
\newacronym[description=Destructive bound enhanced matheuristic]{DBMH}{DBMH}{destructive-bound-enhanced matheuristic}
\newacronym[description=Driver Routing and Scheduling Problem with Multiple Synchronization Constraints]{DRSPMS}{DRSPMS}{Driver Routing and Scheduling Problem with Multiple Synchronization Constraints}

\newacronym[description=no exchange]{noDE}{noDE}{without driver exchange}
\newacronym[description=exchange at regular stops]{DE-RS}{DE-RS}{\underline{d}river \underline{e}xchanges at \underline{r}egular bus \underline{s}tops only}
\newacronym[description=exchange at intermediate stops]{DE-RIS}{DE-RIS}{\underline{d}river \underline{e}xchanges at \underline{r}egular and \underline{i}ntermediate \underline{s}tops}

\newacronym[description=init]{init}{Init}{}
\newacronym[description=init+cLB]{init+cLB}{Init+cLB}{}
\newacronym[description=init+dLB]{init+dLB}{Init+dLB}{}

\newacronym[description=Large Neighborhood Search]{LNS}{LNS}{large neighborhood search}
\newacronym[description=Local Search]{LS}{LS}{local search}
\newacronym[description=Linear Programming]{LP}{LP}{linear programming}

\newacronym[description=Mixed-Integer Program]{MIP}{MIP}{mixed-integer program}

\newacronym[description=Simultaneous Vehicle and Crew Routing and Scheduling Problem]{SVCRSP}{SVCRSP}{simultaneous vehicle and crew routing and scheduling problem}

\newacronym[description=Traveling Salesman Problem]{TSP}{TSP}{traveling salesman problem}

\newacronym[description=Variable MIP Neighborhood Descent]{VMND}{VMND}{variable MIP neighborhood descent}
\newacronym[description=Variable Neighborhood Descent]{VND}{VND}{variable neighborhood descent}
\newacronym[description=Vehicle Routing Problem]{VRP}{VRP}{vehicle routing problem}
\newacronym[description=Vehicle Routing Problem with Multiple Synchronization Constraints]{VRPMS}{VRPMS}{vehicle routing problem with multiple synchronization constraints}
\newacronym[description=Vehicle Routing and Truck Driver Scheduling]{VRTDS}{VRTDS}{vehicle routing and truck driver scheduling}

\newcommand{\graphTSN}{\mathcal{G}}

\newcommand{\metricRemainingWT}{\theta}

\newcommand{\objectiveValue}{z}

\newcommand{\paramContSteering}{T^{cs}}
\newcommand{\paramDailySteering}{T^{ds}}
\newcommand{\paramDailyWorking}{T^{dw}}
\newcommand{\paramBreak}{T^{b}}
\newcommand{\paramConsumption}{c}
\newcommand{\paramNodeTime}{t}
\newcommand{\paramTimeWindow}{tw}
\newcommand{\paramEarliestDep}{e}
\newcommand{\paramLatestDep}{l}
\newcommand{\paramTimeLimit}{3,600}
\newcommand{\paramIntervalLength}{\ell} 
\newcommand{\paramSolutionsKept}{\mu}
\newcommand{\paramSolutionsKeptAbsMin}{\mu^{min}}
\newcommand{\paramNeighborhoods}{\kappa}
\newcommand{\paramDetourLimit}{\zeta}

\newcommand{\paramDeterminism}{p}
\newcommand{\paramTimeLimitGlobal}{\eta}
\newcommand{\paramTimeLimitLBImpr}{\eta^{LB}}
\newcommand{\paramTimeLimitMIPsearch}{\eta^{MIP}}
\newcommand{\paramTimeLimitLocalSearchCB}{\eta^{LS}}
\newcommand{\paramTimeWindowSize}{\vartheta}

\newcommand{\routeDriverCND}{\varrho}
\newcommand{\routeVehicleCND}{\Upsilon}

\newcommand{\solX}{X^{*}}

\newcommand{\setOfRides}{\mathcal{R}}
\newcommand{\setOfDrivers}{\mathcal{K}}
\newcommand{\setOfIntermediateCustNodesPerR}{\mathcal{I}_{r}}

\newcommand{\setOfServiceStationsPerR}{\mathcal{S}_{r}}
\newcommand{\setOfNodesFlatPerR}{\mathcal{N}_{r}}
\newcommand{\setOfNodesFlat}{\mathcal{N}}
\newcommand{\setOfNodesTimeExp}{\mathcal{V}}
\newcommand{\setOfModes}{\mathcal{M}}
\newcommand{\setOfArcsTimeExp}{\mathcal{A}}

\newcommand{\varX}{x}
\newcommand{\varR}{r}
\newcommand{\varB}{b}

\allowdisplaybreaks 

\usepackage{flafter}

\usepackage{dsfont}
\usepackage[scr]{rsfso}

\usepackage{etoolbox}
\patchcmd{\pprintMaketitle}
{\fi\hrule}
{\fi\ifvoid\extrainfobox\else\unvbox\extrainfobox\par\vskip10pt\fi\hrule}
{}{}
\newenvironment{extrainfo}
{\global\setbox\extrainfobox=\vbox\bgroup\parindent=0pt }
{\egroup}
\newsavebox\extrainfobox

\begin{document}

\begin{frontmatter}


\author[tum]{Pia Ammann\corref{cor1}}
\ead{pia.ammann@tum.de}
\ead[url]{http://www.om.wi.tum.de}
\cortext[cor1]{Corresponding author}

\title{Driver Routing and Scheduling with Synchronization Constraints}

\author[tum]{Rainer Kolisch}
\author[tum,mdsi]{Maximilian Schiffer}
\address[tum]{Technical University of Munich, TUM School of Management, Arcisstr. 21, 80333 Munich, Germany}
\address[mdsi]{Munich Data Science Institute, Technical University of Munich, Walther-von-Dyck-Str. 10, 85748 Garching, Germany}

\begin{extrainfo}
	\small
	\textit{Article history:} Received 13 July 2022, Revised 16 December 2022, Revised 19 March 2023, Accepted 26 May 2023.
\end{extrainfo}

\begin{abstract}
This paper investigates a novel type of driver routing and scheduling problem motivated by a practical application in long-distance bus networks. A key difference from other crew scheduling problems is that drivers can be exchanged between buses at arbitrary intermediate stops en route such that our problem requires additional synchronization constraints.
We present a mathematical model for this problem that leverages a time-expanded multi-digraph and derive bounds for the total number of required drivers.
Moreover, we develop a \glsfirst{DBMH} that converges to provably optimal solutions and apply it to a real-world case study for Flixbus, one of Europe's leading coach companies. 
We demonstrate that our matheuristic outperforms a standalone MIP implementation in terms of solution quality and computational time and 
reduces the number of required drivers by up to 56\% compared to approaches currently used in practice.
Our solution approach provides feasible solutions for all instances within seconds and solves instances with up to 390 locations and 70 requests optimally with an average computational time under 210 seconds. 
We further study the impact of driver exchanges on the total driver count and show that allowing for such exchanges leads to average savings of 42.69\%.

\end{abstract}


\end{frontmatter}


\renewcommand\sffamily{}
\section{Introduction}\label{s:intro}
The liberalization of the European long-distance bus market caused a rapid growth in demand for intercity coach travel in recent years, resulting in competition among numerous private coach companies. Despite the long travel duration, intercity bus transportation is generally seen as a low-cost and green alternative to traveling by train or plane. Still, passengers in this sector are susceptible to price and travel duration. 
In this highly competitive environment, a company only survives if it keeps operational costs at a minimum. In this context, driver wages remain a major cost-share. Also, qualified drivers constitute a scarce resource for most bus companies. Accordingly, efficient scheduling and routing that keeps the number of drivers as low as possible remains a crucial planning task.

When creating driver schedules, several legal regulations concerning steering times and rest periods apply. These regulations are particularly important for long-haul passenger transport, where the trip duration often exceeds a driver's maximum steering time without break. For customer satisfaction and to remain competitive, coach companies avoid extended travel times due to drivers taking their legal breaks while passengers are waiting. Instead, they exchange drivers such that a new, sufficiently rested driver continues the trip without any delay. The first driver can then either take a break at the stop where the exchange happened, can return to a depot, or can remain as a passenger on the bus. We refer to the latter as deadheading.
Such driver exchanges can happen at regular bus stops and at any service station along or near the route. Therefore, planning driver schedules in intercity bus networks includes, besides scheduling decisions, additional routing aspects and synchronization constraints. Specifically, the underlying decision problem classifies as a routing and scheduling problem with synchronization constraints and working time regulations.

Against this background, we introduce the \gls{DRSPMS}, which comprises the planning tasks outlined above and is of high practical relevance for operational planning in long-distance bus networks. We develop a \gls{MIP} to formulate the problem on a time-expanded multi-digraph and propose an efficient hybrid solution method that optimally solves real-world instances. 
To separate our work from recent research, we first briefly review related literature in Section~\ref{ss:literatureV2} before further elaborating on our contributions in Section~\ref{ss:contribution}. We then outline the organization of our paper in Section~\ref{ss:orga}.

\subsection{Related Literature}\label{ss:literatureV2}
The \gls{DRSPMS} combines aspects from crew scheduling, truck driver scheduling, and routing with synchronization, which we discuss in the following. In a broader sense, it also relates to routing with intermediate stops for purposes such as replenishment, refueling, or breaks, which has been vividly discussed in recent years. For an overview and classification of this research field, we refer to the comprehensive survey of \cite{SchifferSchneiderEtAl2019}.
%
Related problems in crew scheduling have been thoroughly studied in the past. 
A common assumption in crew scheduling is that trip or flight legs are fixed, such that the main optimization task is to assign crew members to these legs. This assumption does not hold for the \gls{DRSPMS}, in which bus routes are not yet fixed. Instead, one or multiple intermediate stops at arbitrary service stations along or near the route may be inserted between any two regular bus stops. Furthermore, the \gls{DRSPMS} does not require drivers to cover a leg or trip entirely, but allows them to only cover parts of it. Therefore, the \gls{DRSPMS} includes additional routing decisions and synchronization requirements that are not present in crew scheduling problems. Due to these major differences, we do not study the literature on crew scheduling any further. Instead, we refer to \cite{BarnhartCohnEtAl2003} and \cite{DeveciDemirel2018} for a survey on airline crew scheduling, to \cite{HeilHoffmannEtAl2020} for railway crew scheduling, and to \cite{Ibarra-RojasDelgadoEtAl2015} for bus driver scheduling in public transport. 
%
The field of vehicle routing and truck driver scheduling combines scheduling breaks and rest periods for truck drivers with vehicle-related routing decisions but commonly considers vehicles and drivers as inseparable planning units. Contrarily, drivers can switch vehicles at arbitrary intermediate stops en route in the \gls{DRSPMS} and must be treated as separate planning units while still synchronizing their routes. Due to the lack of synchronization requirements in vehicle routing and truck driver scheduling, we do not review work in this research field in more detail but refer to the survey of \cite{KoubaaDhouibEtAl2016} and, for more recent work on this topic, to \cite{GoelIrnich2016}, \cite{SchifferLaporteEtAl2017}, and \cite{TilkGoel2020}.
Table~\ref{tab:literatureOverview_general} compares major characteristics of the \gls{DRSPMS} to related personnel scheduling problems, namely \gls{CS} and \gls{VRTDS}. 
\begin{figure}[]
	\begin{floatrow}
		\ttabbox{%
			\small
			\begin{tabular}{ll|ccc|ccc}
				\toprule							
				&& Assignment & Scheduling & Routing & Working hours & Intermediate stops & Synchronization\\
				\midrule 
				\gls{CS} && \checkmark & \checkmark & & \checkmark & & \\
				\gls{VRTDS} && \checkmark & \checkmark & \checkmark & \checkmark & \checkmark & \\
				\gls{DRSPMS} && \checkmark & \checkmark & \checkmark& \checkmark & \checkmark & \checkmark \\
				\bottomrule	
			\end{tabular}
		}{%
			\caption{Related crew and driver scheduling problems}%
			\label{tab:literatureOverview_general}%
		}
	\end{floatrow}
	
\end{figure}

Due to the synchronization requirements, which result from the possibility of driver exchanges and deadheading, the \gls{DRSPMS} falls into the category of routing problems with synchronization constraints. This relatively new research field has gained increasing attention during the past two decades. 
For a comprehensive review and classification of synchronization in vehicle routing up to 2011, we refer to \cite{Drexl2012}. 
As pointed out by \cite{Drexl2012}, existing publications on vehicle routing with synchronization focus primarily on load synchronization (e.g., for transshipments) and on temporal synchronization of customer visits. Research addressing the temporal and spatial synchronization of non-autonomous vehicles to enable movements is still scarce. 
This so-called movement synchronization is relevant in various settings and applications, including N-echelon routing problems \citep[see][]{GuastarobaSperanzaEtAl2016}, active-passive vehicle routing problems \citep[see][]{MeiselKopfer2014, TilkBianchessiEtAl2018}, and the routing of workers and vehicles \citep[see][]{KimKooEtAl2010, FinkDesaulniersEtAl2019}.
In the following, we survey only literature that considers movement synchronization in the context of vehicle and crew routing and scheduling as contained in the \gls{DRSPMS}. 

In simultaneous vehicle and crew routing and scheduling problems, vehicles and crews are non-autonomous objects who pair up to fulfill a given set of pickup and delivery requests. Furthermore, crews can usually change vehicles at a limited number of locations. Consequently, the routes of crews and vehicles must be synchronized. 
\cite{HollisForbesEtAl2006} study a simultaneous vehicle and crew routing and scheduling problem motivated by an application at Australia Post. 
In their problem setting, postal workers can switch vehicles only at depots and not during a trip. The authors present a set covering formulation and propose a  column generation heuristic to solve the problem of finding minimum cost routes and schedules.
They conduct a computational study with three real-world instances including up to 339 locations and 1,181 shipment requests and show that their heuristic terminates with an average optimality gap of 2.47\% over all instances.
\cite{KergosienLenteEtAl2011} address a joint vehicle and crew routing and scheduling problem arising in medical care. Ambulance drivers, physicians, and vehicles of different types must be routed at minimum costs to transport patients between various care units. Crews can change vehicles at depots only. The authors propose a tabu search algorithm to solve the problem and apply it to instances based on real-world data of a French hospital including 130 transportation requests per day on average. 
Simultaneous vehicle and crew routing and scheduling also arises in long-distance road transport. \cite{DrexlRieckEtAl2013} consider a problem setting in which driver schedules are subject to the European working time regulations, and drivers can only change vehicles at a predefined set of relay stations after completing a daily rest. Furthermore, drivers can use shuttle services to travel between relay stations. 
Routes and schedules have to be determined for trucks and drivers to minimize fixed and variable costs.
The authors develop a two-stage heuristic approach that first fixes truck movements and relay stations and subsequently decides on driver movements. Both stages are solved employing a large neighborhood search.
A case study with real-world data of a major German freight forwarder comprising 2,800~shipment requests between 1,975~customer locations and 157 additional relay stations did not show notable cost savings resulting from the joint planning of vehicle and driver routes. However, the authors state that this result might be due to the unsatisfactory performance of their algorithm and is valid only for their specific setup.
\cite{Dominguez-MartinRodriguez-MartinEtAl2018} study a simultaneous vehicle and crew routing and scheduling problem motivated by a real-world problem of local air traffic in the Canary Islands. Crews can switch aircraft at a limited number of airports, referred to as exchange locations. 
Aircraft and crew routes have to be determined such that distance-dependent costs are minimized.
The authors formally define the problem through a mathematical model with additional valid inequalities and propose a branch-and-cut algorithm to solve the problem. They show in a numerical study on randomly generated test instances that their solution approach can solve instances with up to 30 nodes within a computation time of two hours. 
\cite{LamVanHentenryckEtAl2020} consider a similar problem inspired by an application in military air transportation. Crews can interchange vehicles at different customer locations and may travel as passengers before or after their duty. 
\cite{LamVanHentenryckEtAl2020} formulate an \gls{MIP} model and a constraint programming model to 
find routes and schedules that minimize the number of vehicles and crews as well as the total distance traveled.
They further present two matheuristics that combine both models with a large neighborhood search. 
Computational experiments on randomly generated instances, including up to 11 locations and 50 pickup-and-delivery pairs, show that the constraint-programming-based large neighborhood search outperforms \gls{MIP}-based methods in terms of solution quality. The results also demonstrate that simultaneous routing of crews and vehicles yields notable benefits over sequential approaches. The possibility of vehicle interchanges is essential to obtain good solutions in this specific setting.

As can be seen, simultaneous vehicle and crew routing and scheduling problems arise in various applications. While all existing works consider crew members and vehicles as separate planning units and thus incorporate synchronization constraints, additional constraints vary with the underlying application. 
The possibility of driver exchanges is usually limited to small subsets of locations in order to reduce the problem complexity. Furthermore, many authors consider only simplified or problem-specific working time regulations, such as a maximum route duration.
By contrast, the \gls{DRSPMS} allows driver exchanges at all customer locations and arbitrary intermediate stops, and considers the official working time regulations imposed by the European Union. 
Note that the \gls{DRSPMS} benefits from the fact that vehicle routes are partly fixed as the bus stops and their order are given. Therefore, the routing decisions for vehicles only concern the selection and placement of intermediate stops, which significantly reduces the problem complexity compared to common routing problems. At the same time, however, the large number of potential locations for driver exchanges increases the problem complexity considerably compared to classical crew scheduling problems.
%
To summarize, the \gls{DRSPMS} considered in this paper differs from existing literature in the field of simultaneous vehicle and crew routing and scheduling in several aspects. Table \ref{tab:literatureOverview} compares the characteristics of the \gls{DRSPMS} against the most closely related existing works. Despite its high relevance, no solution algorithm exists that captures all characteristics of the \gls{DRSPMS}, namely i)~driver exchanges at all customer locations, ii)~driver exchanges at intermediate stops, iii)~driver deadheading, and iv)~European driving time regulations. Specifically, we are not aware of any publication considering driver exchanges at both regular and intermediate stops.

\begin{figure}[]
	\begin{floatrow}
		\ttabbox{%
			\small
			\begin{tabular}{lcccccccccccccc}
				
				\toprule							
				&& 
				\parbox{.8cm}{\cite{HollisForbesEtAl2006}}
				&& 
				\parbox[r]{1cm}{\cite{KergosienLenteEtAl2011}} 
				&& 
				\parbox[r]{.8cm}{\cite{DrexlRieckEtAl2013}}
				&& 
				\parbox[r]{1.8cm}{Dominguez-Martin et al. (2018)}
				&& 
				\parbox[r]{.8cm}{\cite{LamVanHentenryckEtAl2020}}
				&& 
				\parbox[r]{.9cm}{This paper\\}
				&
				\\
				\midrule 
				\parbox{3.5cm}{Driver exchanges at \\the depot} && \checkmark && \checkmark &&  &&  &&  && \\
				&&&&&&&&&&&&&\\
				\parbox{3.5cm}{Driver exchanges at \\customer locations} &&  &&  &&  && (\checkmark) && \checkmark && \checkmark\\
				&&&&&&&&&&&&&\\
				\parbox{3.5cm}{Driver exchanges at \\intermediate locations} &&  &&  && (\checkmark) &&  &&  && \checkmark\\
				&&&&&&&&&&&&&\\
				Driver deadheading &&  &&  &&  && \checkmark && \checkmark && \checkmark\\
				&&&&&&&&&&&&&\\
				\parbox{3.5cm}{EU driving time \\regulations} &&  &&  && \checkmark &&  &&  && \checkmark\\
				\bottomrule	
				&&&&&&&&&&&&&\\
				\multicolumn{14}{r}{\footnotesize \textit{Notes.} \checkmark: considered, (\checkmark): partially considered}
				
			\end{tabular}
		}{%
			\caption{Key references on movement synchronization in vehicle and crew routing and scheduling}%
			\label{tab:literatureOverview}%
		}
	\end{floatrow}
\end{figure}

\subsection{Contribution}\label{ss:contribution}
In this paper, we close the research gap outlined in Section \ref{ss:literatureV2} by developing a scalable matheuristic that generates provably optimal solutions for a novel routing and scheduling problem motivated by an application in long-distance bus networks.
Specifically, the contribution of our work is threefold: 

First, we introduce a new type of driver routing and scheduling problem with synchronization constraints. 
To the best of our knowledge, we are the first to consider driver routing and scheduling with synchronization constraints in the context of long-distance bus networks, where, contrary to other applications, drivers can change vehicles not only at customer stops but additionally at arbitrary intermediate stops en route. We develop a mathematical model of our problem based on a time-expanded multi-digraph and introduce valid inequalities to reduce the solution space. We further derive bounds for the total number of required drivers.

Second, we present a hybrid optimization algorithm that, motivated by the requirements of a real-world application, combines mathematical programming techniques with heuristic elements.
Precisely, we alternate between a local search phase to generate solutions and an \gls{MIP} implementation to prove optimality. 
This matheuristic enables us to provide feasible solutions within seconds even for large-scale instances while allowing for statements about the quality of obtained solutions.
We further enrich this approach with a destructive bound improvement phase to tighten lower bounds and speed up convergence. 

Third, we evaluate our solution approach in a comprehensive computational study with real-world data.
The results show that our algorithm is readily applicable in practice. It provides feasible solutions for all instances and solves a huge variety of real-world instances to optimality in 205 seconds on average. 
Moreover, we derive insights for practitioners. 
Compared to current solutions deployed in practice, we can reduce the number of drivers required by up to 56\%. We further show that the consideration of driver exchanges is beneficial in long-distance bus networks and leads to average savings of 43\% per instance.

\subsection{Organization}\label{ss:orga}
The remainder of this paper is structured as follows. First, we describe the problem setting for the \gls{DRSPMS} in Section \ref{s:problemDescription}. In Section \ref{s:solutionMethod}, we introduce a mathematical formulation of the problem and present our matheuristic. Section \ref{s:compStudy} specifies the experimental design for our real-world case study. We present computational results and provide managerial insights in Section \ref{s:results}. Finally, Section \ref{s:conclusion} concludes the paper by summarizing its main insights and suggesting directions for future research.

\section{Problem Description}\label{s:problemDescription}
In the following, we introduce our planning problem. We first describe the practical application before specifying the driving and working time regulations that apply in long-distance bus networks in Europe. Finally, we formally define the \gls{DRSPMS}.

\subsection{Intercity Bus Driver Routing and Scheduling in Practice}\label{ss:ICinPractice}
The deregulation of the European intercity bus market in the past decade has paved the way for new business models in this sector. 
Nowadays, many of the large coach companies in Europe are platform providers that outsource operational tasks to local bus companies, referred to as subcontractors in the following. 
These subcontractors provide buses and drivers and are responsible for executing trips. The platform provider maintains the digital marketplace and takes care of management duties, such as demand forecasting, network planning, bus assignments, or routing and scheduling of buses and drivers.
In this paper, we consider the planning tasks of the platform provider, specifically the planning of driver routes and schedules.

Before creating driver schedules, the platform provider has to plan lines and rides. A line is a sequence of customer stops offered regularly. Once having defined a set of lines to be operated, the platform provider schedules the rides. A ride associates a line with a specific date and time. The platform provider does not fix ride departure times throughout the planning process but instead defines time windows based on target departure times resulting from expected customer demand. Here, we refer to a segment between two subsequent customer stops as a ride segment. 
Given a set of rides and the respective time windows, the platform provider defines driver routes and schedules such that every ride segment is covered by a driver. However, not necessarily the same driver has to cover the entire ride, but drivers may be swapped during the trip. Due to a persistent driver shortage, the platform provider is mainly interested in minimizing the total number of required drivers. Subcontractors often pursue additional objectives such as minimizing the overall workload to ensure driver happiness. Therefore, they might post-optimize driver routes constrained by the total driver count set by the platform provider. As this planning step is not within the scope of the platform provider, we do not consider it any further but refer to Online Appendix~A for a description of how to integrate the subcontractors' objectives into the model developed in this paper.

In many countries, driver schedules are subject to strict driving time regulations that aim to protect driver working conditions and improve road safety.
Such regulations limit, for example, the daily driving time or the maximum driving time without a break.
However, in long-distance bus networks, many lines have a travel duration that exceeds these limits. For coach companies, having the driver take their break while passengers are waiting on the bus is unfavorable as it extends travel times and thus decreases customer satisfaction. Further, short travel times and direct connections with few stops are crucial to gain and sustain market share. Even for very long lines, including breaks can be disadvantageous as most customers travel only for a particular segment on the bus. Customers traveling a short ride segment (of, e.g., 2 hours) might not want their travel time to be interrupted by a break. 
Operating such lines permanently with two drivers could avoid breaks and extended travel times for passengers but would heavily increase operational costs.

As cost-efficient operations are essential for a coach company's success, coach companies prefer to exchange drivers instead of double-booking drivers on buses in order to reduce operational cost. 
Here, a sufficiently rested driver replaces an exhausted driver once she reaches her driving or working time limit. Contrary to other crew scheduling applications, these driver exchanges may not only take place at regular customer stops but additionally at arbitrary service stations along or near the route. Here, driving a short detour to reach a service station is acceptable. Thus, both the drivers' routes and the exact vehicle route for a ride must be defined. 
Specifically, the coach company has to assign drivers to bus rides and schedule the driver's working, driving, and break times. For each ride, potential intermediate stops and exact departure times must be set. Note that we use the terms ride, vehicle route, and bus route interchangeably in the remainder of the paper.
While drivers may use any means of transport to change their location, coach companies favor drivers to travel as passengers on their own bus lines to save costs. Therefore, we assume for the \gls{DRSPMS} that drivers can only move by deadheading or steering a bus of their company.

\subsection{Legal Framework in the European Union}\label{ss:legalFramework}
We now specify the legal regulations that apply in our problem setting.
For countries in the European Union, driving and working times of long-haul bus or truck drivers are regulated by the Regulation (EC) No 561/2006 \citep{EuropeanUnion2006}. 
Specifically, after a maximum steering period of 4.5 hours, drivers have to take a break of at least 45 minutes. This break can be split into two intervals of 15 minutes and 30 minutes. The total steering time per day is limited to 9 hours, whereas the daily working time must not exceed 13 hours. If at least two drivers are present for the entire time between two daily rest periods, the team's working time can be increased up to 21 hours within a 30 hours period. Furthermore, drivers must not drive for more than 56 hours between two consecutive weekly rest periods and no more than 90 hours within two weeks. Additionally, numerous exceptions and application-specific rules apply. 

As we aim to provide a proof of concept with this paper, we assume a planning horizon of a single day for the \gls{DRSPMS} and therefore only consider rules applying to drivers' daily working and steering times. 
Furthermore, we neglect the possibility of splitting breaks. However, integrating these rules into the basic model is a straightforward extension that we leave for future work. 
Table \ref{tab:EUlegislation} gives an overview of the parameters relevant for the \gls{DRSPMS}.

\begin{figure}
	\begin{floatrow}
		\ttabbox{%
			\small
			\begin{tabular}{lll} 					
				\toprule							
				Parameter & Value & Description \\
				\midrule 
				$\paramContSteering$ & 4.5 hours & Maximum steering time without taking a break\\
				$\paramDailySteering$ & 9 hours & Maximum steering time per day (between two consecutive rest periods)\\
				$\paramDailyWorking$ & 13 hours & Maximum working time per day (between two consecutive rest periods)\\
				$\paramBreak$ & 45 minutes & Minimum break duration\\
				\bottomrule	
			\end{tabular}
		}{%
			\caption{Parameters imposed by the European Regulation (EC) No 561/2006}%
			\label{tab:EUlegislation}%
		}
	\end{floatrow}
\end{figure}

\subsection{Problem Definition}\label{ss:modelingConcept}
After describing the practical motivation and the legal framework, we now formally define the \gls{DRSPMS} and summarize all relevant notation in Table \ref{tab:notation}.

\textit{Notation:} Consider a set of rides $\setOfRides$ for a single day. Each ride $r \in \setOfRides$ is associated with a pickup node $p_{r}$ and a delivery node $d_{r}$, which represent the first and the last customer stop of a ride, respectively. Furthermore, we associate a ride $r$ with a set of intermediate customer nodes $\setOfIntermediateCustNodesPerR$ that must be visited in a predefined order. This set is empty for rides consisting only of a single ride segment. In addition to the set of customer stops, defined as $\mathcal{C}_{r} = \{p_{r}\} \: \cup \: \setOfIntermediateCustNodesPerR \: \cup \: \{d_{r}\}$, we consider a set of service stations $\setOfServiceStationsPerR$ for all rides $r \in \setOfRides$ and a central depot $v=0$, where a set of fully rested drivers $\setOfDrivers$ is available. From this virtual depot, drivers can move instantaneously to any other location and back.
Service stations can be used for intermediate stops, for example, to exchange drivers. Here, we consider only service stations that can be reached within a detour limit $\paramDetourLimit \ll \paramBreak$. 
We define node sets~$\setOfNodesFlatPerR$ that include all nodes relevant for ride $r \in \setOfRides$ (i.e., $\setOfNodesFlatPerR = \mathcal{C}_{r}  \cup \setOfServiceStationsPerR$). Furthermore, $\setOfNodesFlat$ denotes the set of relevant nodes over all rides, including the virtual depot (i.e., $\setOfNodesFlat = \cup_{r \in \setOfRides} \setOfNodesFlatPerR \cup \{0\}$).
We define time windows for all nodes~$v \in \setOfNodesFlat \backslash \{0\}$ as $\paramTimeWindow_{v}=[\paramEarliestDep_{v}, \paramLatestDep_{v}]$. These time windows implicitly define the precedence relations between customer nodes of the same ride and are of equal length $\paramTimeWindowSize \ll \paramBreak$. 
Note that the time window size~$\paramTimeWindowSize$ must be greater than or equal to the detour limit~$\paramDetourLimit$. Otherwise, arrival times might be infeasible. To keep deviations from the minimum ride duration at a reasonable level, time window size~$\paramTimeWindowSize$ and detour limit~$\paramDetourLimit$ have to take values significantly smaller than the minimum break duration~$\paramBreak$.

\textit{Solution and Objective Function:} We aim to route and schedule drivers such that the total number of drivers is minimized while serving all ride requests and complying with relevant legal regulations. 
Specifically, a solution i)~determines intermediate stops at service stations, ii)~sets exact departure times at all nodes, iii)~assigns drivers to ride segments, iv)~specifies driver exchanges and deadheading, and v)~sets breaks for all drivers.

\textit{Constraints:} A solution to the \gls{DRSPMS} is subject to the following constraints:
\begin{compactitem}
	\item Every ride $r \in \setOfRides$ is covered by one or multiple drivers such that every customer node $v \in \mathcal{C}_{r}$ is visited within its time windows $\paramTimeWindow_{v}$.
	\item A driver can only cover a single ride at a time.
	\item Drivers can move only while being paired up with a vehicle and vice versa.
	\item Movements of drivers and vehicles are synchronized in time and space.
	\item Driver schedules comply with the legal regulations, specifically the limits $\paramContSteering$, $\paramDailySteering$, and $\paramDailyWorking$ are not exceeded and the minimum break duration $\paramBreak$ is fulfilled.
	\item Detours via service stations must not exceed the detour limit $\paramDetourLimit$.
\end{compactitem}

\section{Methodology}\label{s:solutionMethod}
In the following, we present a methodological framework to solve the \gls{DRSPMS}.
First, we introduce a graph representation in Section \ref{ss:graphRepr}, which allows us to formulate a compact \gls{MIP} in Section \ref{ss:mip}. 
We further enhance this \gls{MIP} by adding valid inequalities and computing lower and upper bounds. We elaborate the bound computations in Section \ref{ss:bounds}. Finally, we describe our matheuristic, which combines local search with mathematical programming techniques and a destructive bound improvement procedure in Section \ref{ss:solApproach}.

\subsection{Time-Expanded Directed Multi-Graph}\label{ss:graphRepr}
In this section, we develop a graph representation that allows us to integrate some of the problem attributes, such as precedence constraints, the possibility of deadheading, or time windows, directly into the graph itself and thereby capture part of the problem complexity outside of the model.

We consider a time-expanded directed multi-graph $\graphTSN = (\setOfNodesTimeExp, \setOfArcsTimeExp)$ with a set of time-expanded nodes  $\setOfNodesTimeExp$ and a set of arcs $\setOfArcsTimeExp$. Every node $i \in \setOfNodesTimeExp$ is associated with a non-time-expanded node $v \in \setOfNodesFlat$ and a discrete time~$\paramNodeTime_{i}$. 
We derive relevant points in time based on the interval length $\paramIntervalLength$, which determines the granularity of time discretization in our graph and thus the approximation's accuracy to the time-continuous problem. 
While large intervals keep the number of time-discrete points low, thus reducing the problem complexity and model quality, shorter time intervals allow for a high-quality approximation but increase the computational effort~\citep[cf.][]{BolandHewittEtAl2019}. 
To overcome this trade-off, \cite{BolandHewittEtAl2017} introduced the concept of dynamic discretization discovery for the continuous-time service network design problem. By iteratively discovering new points in time and adding them to a partially time-expanded network, the algorithm allows solving time-continuous problems to optimality by means of a time-discrete network \citep[cf.][]{BolandHewittEtAl2017, BolandSavelsbergh2019}. Applying dynamic discretization discovery to solve the \gls{DRSPMS} is an interesting extension that we leave for future work.

For the \gls{DRSPMS}, we choose the interval length~$\paramIntervalLength$ such that the time window size~$\paramTimeWindowSize$ is a multiple of it. 
We then create $\frac{\paramTimeWindowSize}{\ell} + 1$ copies of every customer node~$v \in \cup_{r \in \setOfRides} \mathcal{C}_{r}$. 
For service stations, we compute relevant times based on the preceding customer node and create a maximum of $\frac{\paramTimeWindowSize}{\ell}$ copies per node. To obtain the subset of time-expanded nodes associated with a non-time-expanded node~$v \in \setOfNodesFlat$, we define a mapping function~$\mathscr{T}(v): \setOfNodesFlat \rightarrow $ 
$\mathbb{P}(\mathcal{V})$. Here, $\mathbb{P}(\mathcal{V})$ denotes the power set of $\setOfNodesTimeExp$, including all subsets of $\setOfNodesTimeExp$ with cardinality less than or equal to $\frac{\paramTimeWindowSize}{\ell} + 1$.
To represent all feasible connections between time-expanded nodes~$i, j \in \setOfNodesTimeExp$, we define the arc set $\setOfArcsTimeExp$, which includes depot arcs, steering arcs, deadheading arcs, and waiting arcs.
Depot arcs include connections from and to the virtual depot~$0$ for every node $i \in \setOfNodesTimeExp \setminus \{0\}$, whereas steering arcs and deadheading arcs represent connections on which drivers might be steering or deadheading, respectively. Waiting arcs connect nodes associated with the same physical location. 
As there might be multiple arcs between any pair of nodes~$i, j \in \setOfNodesTimeExp$, one for steering and one for deadheading, we introduce a set of modes~$m \in \setOfModes$ with $\setOfModes = \{0, 1\}$ to differentiate arc types. We assign mode~$m=1$ to all steering arcs and mode~$m=0$ to all arcs on which drivers are not steering, that is, depot arcs, waiting arcs, and deadheading arcs. Thus, we characterize an arc by a triple~$(i, j, m)$ with origin~$i \in \setOfNodesTimeExp$, destination~$j \in \setOfNodesTimeExp$, and mode~$m \in \setOfModes$.
As described in Section~\ref{ss:legalFramework}, every driver has a certain amount of continuous steering time available, limited to $T^{cs}$. Steering reduces the remaining amount of this resource by the respective driving time. By taking a break of at least $T^{b}$, drivers can renew their steering time resource. Hence, we associate every arc~$(i,j,m) \in \setOfArcsTimeExp$ with a resource consumption~$\paramConsumption_{ijm}$. 
For steering arcs, the resource consumption is equivalent to the corresponding driving time. For all other arcs, it is equal to zero or may even be negative for arcs which are sufficient to take a break and thus renew the steering time resource. 
We further use the set of arcs~$\setOfArcsTimeExp$ to implicitly model precedence constraints between customer nodes, the detour limit~$\paramDetourLimit$, and the maximum continuous steering time~$T^{cs}$. Accordingly, an arc~$(i,j,m) \in \setOfArcsTimeExp$ only connects two vertices~$i, j \in \setOfNodesTimeExp$ if one of the following conditions holds: i)~exactly one of the nodes represents the depot, ii)~both nodes are customer nodes and node~$j$ is the immediate successor of $i$ , iii)~node~$j$ is a service station relevant for the segment between customer node $i$ and its immediate successor, iv)~node $i$ is a service station relevant for the segment between customer node~$j$ and its immediate predecessor. Here, we consider only service stations that comply with the detour limit~$\paramDetourLimit$. Further, we include without loss of generality only arcs with a travel duration less than or equal to $T^{cs}$ in our arc set~$\setOfArcsTimeExp$.

\begin{figure}
	\begin{floatrow}
		\ffigbox[.3\textwidth]{%
			\centering
			\small
			\begin{tikzpicture}
			
			\tikzstyle{vertex} = [circle, draw=black]
			\tikzstyle{edge} = [->, thick]
			\tikzstyle{selected vertex} = [vertex, fill=red!50]
			
			\node[vertex] (i) at (0,0) {$i$};
			\node[vertex] (j) at (3,3) {$j$};
			\node[vertex] (s) at (2.5,1) {$s$};
			\node[rectangle] (text) at (0,-0.7) {$\paramTimeWindow_{i}=[e_{i}, l_{i}]$};
			\node[rectangle] (text) at (3,3.7) {$\paramTimeWindow_{j}=[e_{j}, l_{j}]$};
			
			\draw[->] (i)--(j);
			\draw[->] (i)--(s) ;
			\draw[->] (s)--(j);
			\end{tikzpicture}
		}{%
			\caption{Ride example}%
			\label{fig:multiGraph}
		}
		\ffigbox[.7\textwidth]{%
			\centering
			\small
			\begin{tikzpicture}
			
			\tikzstyle{vertex} = [circle, draw=black]
			\tikzstyle{edge} = [->, thick]
			\tikzstyle{selected vertex} = [vertex, fill=red!50]
			
			\node[vertex] (i) at (-3,0) {$i_1$};
			\node[vertex] (j) at (0,3) {$j_1$};
			\node[vertex] (s) at (-0.5,1) {$s_1$};
			\node[vertex] (i2) at (3,0) {$i_2$};
			\node[vertex] (j2) at (6,3) {$j_2$};

			\draw[->] (i)--(j) node[midway, left]{$c_{i_1j_11}$};
			\draw[->] (i)-- node[anchor=north]{$c_{i_1s_11}$} (s);
			\draw[->] (s)-- node[pos=.25, above left]{$c_{s_1j_21}$} (j2);

			\path [->, black, densely dotted, bend left=50] (i) edge node[midway, above left]{$c_{i_1j_10}$} (j) ;
			\path [->, black, densely dotted, bend right=60] (i) edge node[anchor=north]{$c_{i_1s_10}$} (s); 
			\path [->, black, densely dotted, bend right] (s) edge node[anchor=south]{$c_{s_1j_20}$} (j2.west);
			
			\draw[->] (i2)--(j2) node[midway, right]{$c_{i_2j_21}$};
			
			\path [->, black, densely dotted, bend right=50] (i2) edge node[anchor=north west]{$c_{i_2j_20}$} (j2);
			\path [->, black, dashed] (i) edge node[pos=.65, below right]{$c_{i_1i_20}$} (i2);
			\path [->, black, dashed] (j) edge node[anchor=south]{$c_{j_1j_20}$} (j2);
			
			\end{tikzpicture}
			\vspace{.2cm}
		}{%
			\caption{Graph representation}%
			\label{fig:TSN}
		}
	\end{floatrow}
\end{figure}

We illustrate our graph representation in a simplified example. Consider customer nodes~$i$ and $j$ associated with time windows $\paramTimeWindow_{i}$ and $\paramTimeWindow_{j}$, respectively, and a ride segment from customer node~$i$ to customer node~$j$. Additionally, there is a service station node $s$. As depicted in Figure \ref{fig:multiGraph}, there are two possibilities to serve this ride segment, either driving directly from $i$ to $j$ or driving via service station $s$. 
Figure \ref{fig:TSN} shows the time-expanded multi-digraph for our example. For the sake of simplicity, we assume $\paramIntervalLength = \paramTimeWindowSize$, that is, we create only two copies of customer node $i$, for times $\paramNodeTime_{i} \in \{\paramEarliestDep_{i}, \paramLatestDep_{i}\}$ such that $\mathscr{T}(i) = \{i_1, i_2\}$. The same holds for node $j$, respectively, such that $\mathscr{T}(j) = \{j_1, j_2\}$.
As nodes $i_1, i_2 \in \setOfNodesTimeExp$ refer to the same physical location, we add a waiting arc from $i_1$ to $i_2$. For the same reason, we add a waiting arc from $j_1$ to $j_2$.
To determine the time-copies of service station $s$, we consider its predecessor, customer node $i$. Given the set of relevant times for node $i$, we can compute the relevant times for $s$ by simply adding the driving time from node $i$ to node $s$. In general, driving via a service station leads to a detour. Therefore, we can never connect the first time-copy of a service station with the first time-copy of its succeeding customer node. Thus, the number of time-copies of a service station is strictly smaller than the number of time-copies needed for a customer node. 
Only a single copy of the service station is needed in this example. Departing from node $i$ at the earliest time possible and driving via service station $s$, we cannot reach node $j$ at time $\paramEarliestDep_{j}$. The next time copy that we can reach is node $j_2$. Thus, we connect node $s_1$ with node $j_2$. When starting at node $i_2$, driving via a service station is not feasible as node $j_2$ cannot be reached. Hence we do not create a second copy of node~$s$. 
Note that we also connect node $s_1$ with node $j_2$ if the detour is smaller than the interval length $\paramIntervalLength$, i.e., we round up the actual arrival time to the next time-discrete value in our graph.
Finally, we add all remaining feasible steering and deadheading arcs, represented by solid and dotted lines, respectively. For the sake of simplicity, we do not display depot arcs.

\subsection{Mathematical Model}\label{ss:mip}
Using the graph representation developed in the previous section, we can now formulate a compact \gls{MIP} model for the \gls{DRSPMS}. We first introduce a basic formulation before we present additional inequalities to break symmetries and tighten the model.

We define the \gls{MIP} on the time-expanded multi-digraph $\graphTSN=(\setOfNodesTimeExp, \setOfArcsTimeExp)$ described in Section \ref{ss:graphRepr}.
To denote out- and ingoing arcs from and into all nodes in a subset $\mathscr{S} \subset \setOfNodesTimeExp$, we define cut sets $\delta^{+}(\mathscr{S}) = \{(i, j, m) \in \setOfArcsTimeExp: i \in \mathscr{S}, j \not\in \mathscr{S}, m \in \setOfModes\}$ and $\delta^{-}(\mathscr{S}) = \{(j, i, m) \in \setOfArcsTimeExp: i \in \mathscr{S}, j \not\in \mathscr{S}, m \in \setOfModes\}$, respectively. Furthermore, we define cut sets for out- and ingoing steering arcs as  $\delta_{1}^{+}(\mathscr{S}) = \{(i, j, m) \in \setOfArcsTimeExp: i \in \mathscr{S}, j \not\in \mathscr{S}, m = 1\}$ and $\delta_{1}^{-}(\mathscr{S}) = \{(j, i, m) \in \setOfArcsTimeExp: i \in \mathscr{S}, j \not\in \mathscr{S}, m = 1\}$. With a slight abuse of notation, we refer to cut sets of singletons $\mathscr{S} = \{i\}$ by $\delta(i)$ instead of $\delta(\{i\})$.
We use binary variables~$\varX_{kijm}$ to state whether driver~$k \in \setOfDrivers$ traverses an arc~$(i,j,m) \in \setOfArcsTimeExp$ ($\varX_{kijm}=1$) or not ($\varX_{kijm}=0$) and derive exact bus routes from the driver routes in a postprocessing step. To track the drivers' steering resource levels, we introduce continuous variables~$\varR_{ki}$, which determine the remaining amount of steering time driver~$k \in \setOfDrivers$ has at node~$i \in \setOfNodesTimeExp$.
With this notation as summarized in Table \ref{tab:notation}, a basic \gls{MIP} for the \gls{DRSPMS} holds as follows:
\begin{figure}
	\begin{floatrow}
		\ttabbox{%
			\small
			\begin{tabularx}{\textwidth}{lX}
				\toprule
				
				\textbf{Sets and Parameters} & \\
				\midrule
				$\setOfRides$ & Set of rides\\
				$\setOfDrivers$ & Set of drivers\\
				
				$\setOfIntermediateCustNodesPerR$ & Set of intermediate customer nodes of ride $r \in \setOfRides$\\
				$\mathcal{C}_{r}$ & Set of customer nodes of ride $r \in \setOfRides, \; \mathcal{C}_{r} = \{p_{r}\} \cup \setOfIntermediateCustNodesPerR \cup \{d_{r}\}$\\
				$\setOfServiceStationsPerR$ & Set of service stations relevant for ride $r \in \setOfRides$ \\
				$\setOfNodesFlatPerR$ & Set of non-time-expanded nodes relevant for ride $r \in \setOfRides$, $\setOfNodesFlatPerR = \setOfServiceStationsPerR \cup \mathcal{C}_{r}$\\
				
				$\setOfNodesTimeExp$ & Set of time-expanded nodes\\
				$\setOfModes$ & Set of modes $m \in \setOfModes$, with $m=1$ for steering arcs and $m=0$ for all other arcs\\	
				$\setOfArcsTimeExp$ & Set of feasible arcs $(i,j,m)$ between time-expanded nodes $i, j \in \setOfNodesTimeExp$ with mode $m \in \setOfModes$\\
				$\mathscr{T}(v)$ & Function mapping non-time-expanded node $v \in \setOfNodesFlat$ to the respective subset of time-expanded nodes, $\mathscr{T}(v): \setOfNodesFlat \rightarrow \mathbb{P}(\mathcal{V})$\\ 
				$\delta^{+}(\mathscr{S}), \delta^{-}(\mathscr{S})$ & Cut sets of outgoing / ingoing arcs from / into subset $\mathscr{S} \subset \setOfNodesTimeExp$\\
				$\delta_{1}^{+}(\mathscr{S}), \delta_{1}^{-}(\mathscr{S})$ & Cut sets of outgoing / ingoing steering arcs from / into subset $\mathscr{S} \subset \setOfNodesTimeExp$\\
				$p_{r}, d_{r}$ & Pickup and delivery nodes of ride $r \in \setOfRides$\\	
				$\paramConsumption_{ijm}$ & Resource consumption associated with arc $(i, j, m) \in \setOfArcsTimeExp$\\
				$\paramNodeTime_i$ & Time associated with node $i \in \setOfNodesTimeExp$\\
				\midrule
				\textbf{Decision Variables} & \\
				\midrule
				$\varX_{kijm}$ & Binary variables that indicate if driver $k \in \setOfDrivers$ traverses arc $(i,j,m) \in \setOfArcsTimeExp$ ($\varX_{kijm}=1$) or not ($\varX_{kijm}=0$)\\
				$\varR_{ki}$ & Continuous variables that track the remaining amount of steering time for driver $k \in \setOfDrivers$ at node $i \in \setOfNodesTimeExp$\\
				\bottomrule
			\end{tabularx}
		}{%
			\caption{Notation}%
			\label{tab:notation}%
		}
	\end{floatrow}
\end{figure}

\begin{equation}
    \textbf{min} \qquad \objectiveValue = \sum_{k \in \setOfDrivers}\sum_{(i,j,m) \in \delta^{+}(\mathscr{T}(0))} \varX_{kijm}\label{eq:objective}
\end{equation}

subject to

\begin{align}
	\sum_{(i,j,m) \in \delta^{+}(\mathscr{T}(0))} \varX_{kijm} \leq 1 \qquad & \forall \; k \in \setOfDrivers\label{eq:driverStart}\\
	\sum_{(j,i,m) \in \delta^{-}(i)} \varX_{kjim} = \sum_{(i,j,m) \in \delta^{+}(i)} \varX_{kijm} \qquad & \forall \; k \in \setOfDrivers,\; i \in \setOfNodesTimeExp\label{eq:driverFlow}\\
	\sum_{i \in \mathscr{T}(v)} \; \sum_{k \in \setOfDrivers} \; \sum_{(i,j,m) \in \delta_{1}^{+}(i)} \varX_{kijm} = 1 \qquad & \forall \; r \in \setOfRides,\; v = p_{r}\label{eq:requestFlowP}\\
	\sum_{i \in \mathscr{T}(v)} \; \sum_{k \in \setOfDrivers} \; \sum_{(j,i,m) \in \delta_{1}^{-}(i)} \varX_{kjim} = 1 \qquad & \forall \; r \in \setOfRides, \; v = d_{r}\label{eq:requestFlowA}\\
	\sum_{k \in \setOfDrivers} \; \sum_{(i,j,m) \in \delta_{1}^{+}(i)} \varX_{kijm} = \sum_{k \in \setOfDrivers} \; \sum_{(j,i,m) \in \delta_{1}^{-}(i)} \varX_{kjim} \qquad &  \forall \; r \in \setOfRides,\; v \in \setOfIntermediateCustNodesPerR \cup \setOfServiceStationsPerR,\; i \in \mathscr{T}(v)\label{eq:requestFlowI}\\
	\varX_{lij0} \leq \sum_{k \in \setOfDrivers \backslash \{l\}} \varX_{kijm} \qquad & \forall \; l \in \setOfDrivers,\; (i,j,m) \in \setOfArcsTimeExp: m=1\label{eq:deadheading}
\end{align}

\vspace{-0.8cm}
\begin{align}
	\sum_{\substack{(i,j,m) \in \setOfArcsTimeExp:\\ m=1}} \paramConsumption_{ijm} \varX_{kijm} \leq \paramDailySteering \qquad & \forall \; k \in \setOfDrivers\label{eq:HOSdailySteering}\\
	\sum_{i \in \mathscr{T}(0)} \sum_{(j,i,m) \in \delta^{-}(i)} \paramNodeTime_{j} \varX_{kjim} - \sum_{i \in \mathscr{T}(0)} \sum_{(i,j,m) \in \delta^{+}(i)}   \paramNodeTime_{j}\varX_{kijm} \leq \paramDailyWorking \qquad & \forall \; k \in \setOfDrivers\label{eq:HOSdailyWorking}\\
	\varR_{kj} \; \leq \; \varR_{ki} - \sum_{\substack{m \in \setOfModes:\\(i, j, m) \in \setOfArcsTimeExp}}  (\paramConsumption_{ijm}  \varX_{kijm}) + \paramContSteering  (1- \sum_{\substack{m \in \setOfModes:\\(i, j, m) \in \setOfArcsTimeExp}} \varX_{kijm}) \qquad & \forall \; k \in \setOfDrivers, \; (i,j): (i,j,m) \in \setOfArcsTimeExp \label{eq:HOSreduceRwhenSteering}
\end{align}

\vspace{-0.8cm}
\begin{align}
	\varX_{kijm} \in \{0,1\} \qquad & \forall \; k \in \setOfDrivers, \; (i,j,m) \in \setOfArcsTimeExp \label{eq:domainX}\\
	0 \leq \varR_{ki} \leq \paramContSteering \qquad & \forall \; k \in \setOfDrivers, \; i \in \setOfNodesTimeExp\label{eq:domainR}
\end{align}

The objective function (\ref{eq:objective}) minimizes the number of drivers leaving the depot, that is, the total number of drivers required. 
Constraint sets~(\ref{eq:driverStart}) - (\ref{eq:requestFlowI}) regulate driver and requests flows in our network. Constraints~(\ref{eq:driverStart}) state that drivers can leave the depot only once and accordingly execute only a single route. Constraints~(\ref{eq:driverFlow}) ensure flow consistency through the network for every driver. Constraint sets~(\ref{eq:requestFlowP}) - (\ref{eq:requestFlowI}) require all requests to be served. While constraint sets~(\ref{eq:requestFlowP}) and (\ref{eq:requestFlowA}) ensure that pickup and delivery nodes of every request~$r \in \setOfRides$ are visited, Constraints~(\ref{eq:requestFlowI}) guarantee the consistency of request flows at intermediate customer stops and at service station nodes. 
Constraints~(\ref{eq:deadheading}) allow deadheading for a driver $l \in \setOfDrivers$ only if there is another driver~$k \in \setOfDrivers \backslash \{l\}$ steering on the same arc. 
To integrate working and driving time regulations, we further define Constraints~(\ref{eq:HOSdailySteering}) - (\ref{eq:HOSreduceRwhenSteering}). Constraints~(\ref{eq:HOSdailySteering}) and (\ref{eq:HOSdailyWorking}) limit the daily steering time and daily working time for each driver. Constraints~(\ref{eq:HOSreduceRwhenSteering}) update the driver's steering resource level after traversing an arc. 
Finally, constraint sets~(\ref{eq:domainX}) - (\ref{eq:domainR}) define the variables' domains.

To improve the \gls{MIP}'s performance, we formulate additional Constraints~(\ref{eq:HOSinitR}) - (\ref{eq:symmetryBreakingT}), which cut off symmetric solutions and tighten some variable assignments. 
\begin{align}
	\varR_{ki} \geq \paramContSteering  (\sum_{(i,j,m) \in \delta^{+}(i)} \varX_{kijm}) \qquad & \forall \; k \in \setOfDrivers, \; i \in \mathscr{T}(0)\label{eq:HOSinitR}\\
	\varR_{ki} \leq \paramContSteering  (\sum_{(i,j,m) \in \delta^{+}(i)} \varX_{kijm}) \qquad & \forall \; k \in \setOfDrivers, \; i \in \setOfNodesTimeExp\label{eq:tighterDomainR}\\
	\sum_{i \in \mathscr{T}(0)} \sum_{(i,j,m) \in \delta^{+}(i)} \varX_{kijm} \geq \sum_{i \in \mathscr{T}(0)} \sum_{(i,j,m) \in \delta^{+}(i)} \varX_{lijm} \qquad & \forall \; k \in \{1,..,|\setOfDrivers|-1\},\; l=k+1\label{eq:symmetryBreaking}\\
	\sum_{i \in \mathscr{T}(0)} \sum_{(i,j,m) \in \delta^{+}(i)}  \paramNodeTime_{i} \varX_{kijm} \geq \sum_{i \in \mathscr{T}(0)} \sum_{(i,j,m) \in \delta^{+}(i)}   \paramNodeTime_{i} \varX_{lijm} \qquad & \forall \; k \in \{1,..,|\setOfDrivers|-1\},\; l=k+1\label{eq:symmetryBreakingT}
\end{align}

Constraints~(\ref{eq:HOSinitR}) initialize the resource level of activated drivers, whereas Constraints~(\ref{eq:tighterDomainR}) fix the resource level $\varR_{ki}$ to zero if driver $k \in \setOfDrivers$ never visits node $i \in \setOfNodesTimeExp$.
Next, Constraints~(\ref{eq:symmetryBreaking}) and (\ref{eq:symmetryBreakingT}) break some of the symmetry in the solution space by activating drivers one after another and assigning routes to drivers by descending departure times.

\subsection{Bound Computations}\label{ss:bounds}
To iterate over drivers in our mathematical model, we calculate an artificial upper bound for the number of drivers.
We compute a simple upper bound by considering the drivers' scarcest resource, the continuous steering time $T^{cs}$, assuming that there is no possibility of renewing this resource, that is, we neglect any possibility of taking a break. We thus divide every ride $r \in \mathcal{R}$ into $s_r$ segments of total driving times less than or equal to $T^{cs}$ and assume that a new driver is assigned to each segment. Then, the total number of segments over all rides gives an upper bound for the number of required drivers, see Equation (\ref{eq:UB}). We use the resulting upper bound $UB$ to define the cardinality of the set of drivers $\setOfDrivers = \{1, .., UB\}$. 
\begin{align}
UB \; = \; \sum_{r \in \setOfRides} s_{r} \label{eq:UB}
\end{align}

Preliminary computations showed the \gls{LP} relaxation of the \gls{MIP} defined in Section~\ref{ss:mip} to be rather weak. Therefore, we derive two constructive lower bounds that further strengthen the \gls{MIP} formulation. For a computational analysis of the different lower bounds, we refer to Online Appendix~B. We compute a first constructive lower bound by relaxing some of the legal regulations, particularly the limit on the continuous steering time $T^{cs}$. We divide the minimum total steering time by the maximum daily steering time per driver ($\paramDailySteering$). Here, the minimum total steering time is the sum of driving times over all direct steering arcs between customer locations, such that our first lower bound $LB_1$ results in:
\begin{align}
	LB_1 \; = \; \lceil{(\sum_{(i,j,m) \in \setOfArcsTimeExp: i, j \in \cup_{r \in \setOfRides}\mathcal{C}_{r}, m=1} \paramConsumption_{ijm})\; /T^{ds}} \rceil{}\label{eq:LB_trivial}
\end{align}
This approach assumes a perfect spatial and temporal fit between all rides, allowing drivers to utilize their available daily steering time perfectly. It also assumes that drivers can always drive directly from one customer stop to the next without a detour via any service station. Further, it disregards any conflicts between rides. Therefore, it constitutes a rather weak lower bound, particularly for instances involving many short rides taking place simultaneously as these cannot be executed by the same driver.
Therefore, we derive a second lower bound $LB_2$ by determining the number of parallel rides $r_{t}$ per time interval $t \in \mathcal{T}$. Here, we first obtain time intervals $t \in \mathcal{T}$, representing minimum overlaps resulting from the latest start and earliest arrival times at regular customer stops of all rides. For every time interval $t \in \mathcal{T}$, we then count the number of rides active within this interval. The maximum ride count defines the minimum number of drivers required simultaneously:
\begin{align}
	LB_2\; = \; \max_{t \in \mathcal{T}}(r_{t})\label{eq:LB_r_parallel}
\end{align}
As neither of the two approaches is provably superior over the other, we always compute both lower bounds and take the maximum value as a global lower bound $LB$: 
\begin{align}
	LB \; = \; \max(LB_1,LB_2)\label{eq:LB}
\end{align}
Given lower bound~$LB$, we can further improve our \gls{MIP} by adding constraints (\ref{eq:LB_x}) - (\ref{eq: VI_in}). Note that lower bounds $LB_1$, $LB_2$, and $LB$ are independent of the level of time discretization and are also valid for the original (time-continuous)  problem as described in Section~\ref{s:problemDescription}. 
\begin{align}
	\sum_{k \in \setOfDrivers} \sum_{i \in \mathscr{T}(0)}\sum_{(i,j,m) \in \delta^{+}(i)}\varX_{kijm} \; \geq \; LB\label{eq:LB_x}\\
	\sum_{i \in \mathscr{T}(0)}\sum_{(i,j,m) \in \delta^{+}(i)}\varX_{kijm} \; \geq \; 1 \qquad & \forall \; k \in \{1,..,LB\} \label{eq: VI_out}\\
	\sum_{i \in \mathscr{T}(0)}\sum_{(j,i,m) \in \delta^{-}(i)}\varX_{kjim} \; \geq \; 1 \qquad & \forall \; k \in \{1,..,LB\} \label{eq: VI_in}
\end{align}
Constraint~(\ref{eq:LB_x}) sets the lower bound for the total number of required drivers. Given this lower bound and our symmetry breaking Constraints~(\ref{eq:symmetryBreaking}) - (\ref{eq:symmetryBreakingT}), we force drivers with an index up to the lower bound to be activated, that is, they must have an active arc to and from the depot, as stated in Constraints~(\ref{eq: VI_out}) - (\ref{eq: VI_in}). 

\subsection{Destructive-Bound-Enhanced Matheuristic}\label{ss:solApproach}
Due to the inherent complexity of the \gls{DRSPMS}, the basic \gls{MIP} implementation is not tractable for real-world instances.
Preliminary experiments showed that commercial \gls{MIP} solvers are able to solve small instances but fail to find feasible solutions for real-world-sized instances even after hours of computations.
Accordingly, motivated by our real-world application, we aim at developing a solution approach that is readily applicable in practice. This approach must provide feasible high-quality solutions within a reasonable amount of computational time. In fact, as saving a single driver may lead to significant cost reductions, we aim to find optimal solutions. 
While branch-and-price approaches are state-of-the-art for routing problems, the presence of synchronization requirements causes considerable difficulties with labeling algorithms commonly used in branch-and-price to solve the pricing problem \citep[cf.][]{Drexl2012}. Alternative branch-and-price methods that are tailored to solving such a pricing problem solve only relatively small instances to optimality \citep[cf.][]{TilkDrexlEtAl2019}. Synchronization constraints also lead to difficulties with classical improvement heuristics that assume independence of routes \citep[cf.][]{Drexl2012}. Thus, it remains an open question which methods are most effective for solving routing problems with synchronization constraints. Accordingly, a large variety of different solution methods often tailored to the specific problem at hand exists. To satisfy the requirements from the real-world application of the \gls{DRSPMS}, we propose a matheuristic that provides feasible solutions quickly through its heuristic components while still converging to provably optimal solutions. To speed up convergence, we further enrich our matheuristic with a destructive improvement phase.
Following recent performance improvements of mathematical programming solvers, matheuristics have become increasingly popular for various problems and application areas. For other works that apply similar concepts, we refer the interested reader to, e.g., \cite{LarrainCoelhoEtAl2017, LamVanHentenryckEtAl2020}, and \cite{YinDArianoEtAl2021}.

Algorithm \ref{alg:framework} outlines the general structure of our \glsfirst{DBMH}. First, we use a construction heuristic to generate a feasible start solution. We further improve this solution by applying a local search procedure. 
\begin{algorithm}[b]
    \caption{Destructive-bound-enhanced matheuristic}
    \label{alg:framework}
	\SetKwInOut{Input}{Input}	
	\Input{Constructive lower bound $LB$, time limits $\paramTimeLimitGlobal$, $\paramTimeLimitLBImpr$,$\paramTimeLimitMIPsearch$, $\paramTimeLimitLocalSearchCB$} 
    $S \leftarrow \; ConstructionHeuristic()$\\
    $S \leftarrow \; LocalSearch(S, \paramTimeLimitGlobal)$\\
    \If{$f(S) = LB$}{\Return{$S$}}
    $LB, S \leftarrow \; DestructiveBoundImprovement(LB, S, \paramTimeLimitLBImpr)$\\
    \If{$f(S) = LB$}{\Return{$S$}}
    \While{$\paramTimeLimitGlobal \; not \; reached$}{
		$S,LB \leftarrow \; SolveMIP(LB, S, \paramTimeLimitMIPsearch)$\\
		$S \leftarrow \; LocalSearch(S, \paramTimeLimitLocalSearchCB)$\\         
    }
    \Return{$S$}
\end{algorithm}

If the objective value of the resulting solution $S$ equals the constructive lower bound $LB$, solution $S$ is optimal, and we terminate the search. 
Otherwise, we invoke a destructive bound improvement procedure to improve the lower bound for our problem. If the initial solution $S$ is proven to be optimal during this phase or an optimal solution $S^{*}$ is found, we stop. 
Otherwise, we enter the hybrid optimization phase, which is inspired by the variable \gls{MIP} neighborhood descent of \cite{LarrainCoelhoEtAl2017}. This approach enriches the implementation of an \gls{MIP} by a heuristic component. Specifically, the underlying idea is to alternate between a local search phase to generate improved solutions and the \gls{MIP} to prove optimality. 
Here, we start from an initial solution $S$ and a destructive lower bound $LB$, and solve the \gls{MIP} defined in Sections \ref{ss:mip} and \ref{ss:bounds} until one of the following three conditions holds: i) optimality is proven, ii) we find a new incumbent, or iii) we reach a time limit. 
Note that we first run the \gls{MIP} without providing a start solution and inject the initial solution $S$ only if the \gls{MIP} solver did not obtain any feasible solution after the time limit $\paramTimeLimitMIPsearch$.
In case of proven optimality or timeout, the algorithm terminates. 
By optimality we mean optimality of the time-discrete problem, which is not necessarily optimal for the original (time-continuous) problem.  Solutions, for which optimality is proven by means of the constructive lower bounds, however, are truly optimal, both for the time-discrete and the time-continuous problem.
If we find a new incumbent, we reinvoke the local search based on the new solution $S$ and continue iterating between the local search phase and the \gls{MIP} phase. Note that the \gls{MIP} phase might also update the current lower bound $LB$.
Contrary to \cite{LarrainCoelhoEtAl2017}, who explore neighborhoods by solving constrained subproblems with an \gls{MIP} solver, we implement all heuristic components outside of the \gls{MIP} implementation. This has proven to be computationally more efficient for the \gls{DRSPMS} as \gls{MIP} solvers still suffer from slow convergence for constrained variants of the original problem and show difficulties finding new feasible solutions quickly. Furthermore, we follow a best-improvement strategy whereas \cite{LarrainCoelhoEtAl2017} switch back to the full \gls{MIP} with the first improvement found.

In the remaining sections, we describe the components of our matheuristic in more detail. First, we outline the destructive bound improvement in Section \ref{ss:destructiveImprovement}. Next, we explain the solution representation used in the heuristic components in Section \ref{sss:solRepr} before we describe the construction heuristic in Section \ref{sss:constrHeuristic}. Finally, we specify the local search phase in Section \ref{sss:vndNeighborhoods}. We provide a computational analysis on the impact of all components on the solution quality in Online Appendix~C.

\subsubsection{Destructive Bound Improvement}\label{ss:destructiveImprovement} 
Preliminary computations in the development phase showed that the hybrid optimization approach enriching the \gls{MIP} implementation by a local search component as in \cite{LarrainCoelhoEtAl2017} had difficulties raising the lower bound and therefore converged rather slowly for many instances.
To tighten our lower bounds and accelerate convergence, we introduce a destructive bound improvement phase based on the idea presented in \cite{KleinScholl1999} for solving the resource-constrained project scheduling problem. \cite{KleinScholl1999} differentiate two main classes for computing lower bounds for minimization problems. Constructive bounds are computed using direct methods, such as solving the \gls{LP} relaxation. Destructive bounds, on the other hand, are computed by imposing a bound on the objective function and then trying to contradict the feasibility of this restricted problem. When setting the upper bound to a low value, the search space of feasible solutions is reduced significantly and might even be empty.   
\cite{KleinScholl1999} propose working with the \gls{LP} relaxation of the original problem when applying destructive techniques. Due to the weak \gls{LP} relaxation of our problem, however, we work with the restricted integer problem instead. 

Algorithm \ref{alg:destructiveImpr} outlines the main steps of this procedure.
\begin{algorithm}[]
    \caption{Destructive bound improvement}
    \label{alg:destructiveImpr}
	\SetKwInOut{Input}{Input}	
	\Input{Constructive lower bound $LB$, initial solution $S$, time limit $\paramTimeLimitLBImpr$} 
    \While{$\paramTimeLimitLBImpr$ not reached}{
		\If{$f(S) = LB$}{ 
			\Return{$S$}
		}
		Build $P(LB): \sum_{k \in \setOfDrivers}\sum_{(i,j,m) \in \delta^{+}(\mathscr{T}(0))} \varX_{kijm} \leq LB$\\ 
		Solve $P(LB)$\\
		\If{feasible solution $S$}{
			\Return{$S$}
		}
		\If{$P(LB)$ is infeasible}{
			$LB \leftarrow LB + 1$\\	
		}         
    }
    \Return{$LB$}
\end{algorithm}

Given an initial constructive lower bound $LB$ for our original problem $P$, we first check whether we can prove the initial solution $S$ to be optimal. If this is not the case, that is, $f(S) > LB$, we create a restricted problem $P(LB)$ by adding the constraint $\sum_{k \in \setOfDrivers}\sum_{(i,j,m) \in \delta^{+}(\mathscr{T}(0))} \varX_{kijm} \leq LB$.
We then try to solve this restricted problem $P(LB)$ using a standard solver. If a feasible solution $S$ is found, we terminate and return $S$ as this is an optimal solution for our original problem $P$.
If the restricted problem is proven infeasible, we increment our lower bound $LB \leftarrow LB+1$ and restart the improvement phase. We repeat these steps until we prove optimality or reach a time limit. If we reach the time limit $\paramTimeLimitLBImpr$ without a feasible solution, the current lower bound $LB$ is a valid lower bound for our original problem $P$.
Note that, in contrast to the constructive lower bounds defined in Section~\ref{ss:bounds}, the destructive lower bound depends on the level of time discretization and thus only holds for the time-discrete problem as defined in Section~\ref{ss:mip}.

\subsubsection{Solution Representation and Evaluation}\label{sss:solRepr}
For the heuristic components of our solution approach, we choose a solution representation that considers only a subset of the decision variables defined in Section \ref{ss:mip}. Specifically, we present a solution as a collection of driver routes. Thus, we model only arcs that are traversed in the current solution (i.e., $(i,j,m) \in \setOfArcsTimeExp: \sum_{k \in \setOfDrivers} \varX_{kijm} \geq 1$). Furthermore, we consider only drivers that execute at least one segment of a ride and refer to this subset of active drivers in a solution $S$ as $\bar{\setOfDrivers}(S) \subseteq \setOfDrivers: \sum_{(i,j,m) \in \setOfArcsTimeExp: m=1} \varX_{kijm} \geq 1 \; \forall \; k \in \bar{\setOfDrivers}$.
We assign every driver $k \in \bar{\setOfDrivers}$ a route $\routeDriverCND_k$ through our time-expanded multi-digraph $\mathcal{G}$. A route consists of an ordered list of arcs $(i, j, m)$ to be traversed by the respective driver, such that $x_{kijm}=1 \; \forall \; (i,j,m) \in \routeDriverCND_k$. 
Departure times $t(S)$ can be derived from the times $t_i$ associated with the nodes $i \in \setOfNodesTimeExp$ that are included in the driver routes of solution $S$. Note that we do not explicitly model the resource variables $\varR_{ki}$ but instead check a solution for resource feasibility after every modification of the driver routes. 
Then, we can derive the objective function value $f(S)$ associated with a solution $S$ from the cardinality of the subset $\bar{\setOfDrivers}(S)$, such that $f(S) = |\bar{\setOfDrivers}|$.
Considering only arcs with mode $m=1$, we can derive the 
exact routes for each bus ride, which we also refer to as vehicle routes~$\routeVehicleCND(S)$.

\subsubsection{Construction Heuristic}\label{sss:constrHeuristic}
To construct an initial solution, we use a greedy two-step approach that first fixes vehicle routes and schedules and, second, defines driver routes and schedules. As the number of drivers is not limited, we can always generate a feasible solution. 
To create vehicle routes, we assume that all ride segments are served without intermediate stops at service stations. We further fix the departure times of all ride segments to the earliest time possible. Based on the resulting vehicle routes, we assign drivers as follows. For every vehicle route, starting from the one with the earliest start time, we apply the following steps:
\begin{enumerate}
	\item If a driver is available at the start node and can execute the first segment, we assign her to the current vehicle route. Otherwise, we add a new driver. We let the selected driver steer until i) the vehicle route is completed or ii) the driver reaches her continuous steering time limit $\paramContSteering$.
	\item In case of i), we let the current driver and all deadheading drivers wait at the terminal node until they can take over another vehicle route or have exhausted their total working time.
	\item In case of ii), we decide whether to let the driver wait at the current stop or whether to stay on the bus (deadheading). If another bus later visits this stop and the waiting time is sufficient to take a break of $\paramBreak$ but no longer than a predefined upper bound, e.g., $4 \paramBreak$, the driver waits. Otherwise, the driver deadheads.
\end{enumerate}
We continue Steps 1 - 3 until we assigned all ride segments. The resulting vehicle routes~$\routeVehicleCND(S)$, departure times~$t(S)$, and driver routes~$\routeDriverCND(S)$ define our start solution~$S$.

\subsubsection{Local Search in a Composite Neighborhood}\label{sss:vndNeighborhoods}
In contrast to \cite{LarrainCoelhoEtAl2017}, who complement their \gls{MIP} with a variable neighborhood descent, we apply a local search based on a composite neighborhood, which has proven to be more efficient for the \gls{DRSPMS}. We apply this local search procedure to improve solutions provided by the construction heuristic and by the \gls{MIP} solver. Composite neighborhoods have already been used successfully by other authors \citep[see, e.g.,][]{SchifferWalther2018, SchneiderLoffler2019}. For computational analysis of the different neighborhood search variants, we refer to Online Appendix~D. 

We define our composite neighborhood through $\paramNeighborhoods$ operators. 
In every iteration of our local search, we select the best solution $S$ in our composite neighborhood.
We continue this search until we reach a time limit or no further improvement is possible, indicating that we found a local optimum with respect to our composite neighborhood, and return the best solution found so far. 

All $\paramNeighborhoods$ operators in our composite neighborhood follow a best-improvement strategy, that is, they modify the current solution until no further improvement is possible. 
In every iteration, multiple new solutions are generated. Any of these new solutions can then be further explored and serves as a basis for creating  multiple modified solutions in the next iteration. However, exploring all these possible modifications would result in an unreasonable computational effort. Moreover, many of the solutions show only minor structural differences but have the same objective function value because small changes, such as shifting the departure time of a single ride segment, do not necessarily affect the total number of drivers needed.
Therefore, we focus only on a small subset of promising solutions. We introduce a secondary objective as a tiebreaker that allows us to differentiate and rank solutions with the same primary objective value. This secondary objective accounts for the remaining amount of working time $\metricRemainingWT(S)$ of drivers who are active in the current solution $S$, considering only the remaining working time before starting or after finishing their route:
\begin{align}
\metricRemainingWT(S) = \sum_{k \in \setOfDrivers} \left( \paramDailyWorking - \left( \sum_{(j,i,m) \in \delta^{-}(\mathscr{T}(0))} \paramNodeTime_{i} \varX_{kjim} - \sum_{(i,j,m) \in \delta^{+}(\mathscr{T}(0))} \paramNodeTime_{i} \varX_{kijm} \right) \right) \label{eq:S_eff}
\end{align} 
Solutions with a high remaining working time $\metricRemainingWT(S)$ might offer a higher potential for consolidating driver routes and thus allow for reducing the number of drivers. 
Hence, we order solutions with identical primary objective function values by decreasing amount of remaining working time and consider only the first $\paramSolutionsKept$(\%) out of this list in each iteration. 
At minimum, we keep an absolute number of $\paramSolutionsKeptAbsMin$ solutions per iteration.

Due to the synchronization requirements present in the \gls{DRSPMS}, standard local search operators are not applicable. Instead, we developed seven problem-specific operators. These operators address different problem attributes, precisely the composition of driver routes $\routeDriverCND(S)$, the composition of vehicle routes $\routeVehicleCND(S)$, and the departure times $t(S)$ in a solution $S$. Table \ref{tab:neighborhoods} provides an overview of these operators, which are defined as follows:

\begin{figure}  
	\begin{floatrow}
		\ttabbox{%
			\small
			\begin{tabular}{lrr}
				\toprule
				Operator & \multicolumn{1}{c}{Modifies} & \multicolumn{1}{c}{Size}\\
				\midrule
				Reassign Segments & $\routeDriverCND(S)$ & $\paramSolutionsKept \sum_{n=1}^{|\mathcal{RS}|-1} \frac{|\mathcal{RS}|!}{(|\mathcal{RS}|-n)!}$ \\[4pt] 
				Postpone & $t(S)$, $\routeDriverCND(S)$ & $\paramSolutionsKept \; 2^{\frac{\paramTimeWindowSize}{\paramIntervalLength} |\mathcal{R}|}$\\[4pt] 
				Prepone & $t(S)$, $\routeDriverCND(S)$ & $\paramSolutionsKept \; 2^{\frac{\paramTimeWindowSize}{\paramIntervalLength} |\mathcal{R}|}$\\[4pt] 
				Insert Intermediate Stop (Random) & $\routeVehicleCND(S)$, $t(S)$, $\routeDriverCND(S)$ & $\paramSolutionsKept \; 2^{\frac{\paramTimeWindowSize}{\paramIntervalLength} |\mathcal{R}|}$\\[4pt] 
				Insert Intermediate Stop (Shortest Detour) & $\routeVehicleCND(S)$, $t(S)$, $\routeDriverCND(S)$ & $\paramSolutionsKept \; 2^{\frac{\paramTimeWindowSize}{\paramIntervalLength} |\mathcal{R}|}$\\[4pt] 
				Insert Intermediate Stop (Synchronization Potential) & $\routeVehicleCND(S)$, $t(S)$, $\routeDriverCND(S)$ & $\paramSolutionsKept \; 2^{\frac{\paramTimeWindowSize}{\paramIntervalLength} |\mathcal{R}|}$\\[4pt] 
				Remove Intermediate Stop & $\routeVehicleCND(S)$, $t(S)$, $\routeDriverCND(S)$ & $\paramSolutionsKept \; 2^{\frac{\paramTimeWindowSize}{\paramIntervalLength} |\mathcal{R}|}$\\[4pt] 
				\bottomrule
				\multicolumn{3}{r}{\footnotesize \textit{Notes.} $\; \mathcal{RS}$: Set of ride segments.} 
			\end{tabular}
		}{%
			\caption{Local search operators}%
			\label{tab:neighborhoods}%
		}
	\end{floatrow}
\end{figure}

\paragraph{Reassign Segments}
This operator assumes fixed vehicle routes $\routeVehicleCND(S)$ and departure times $t(S)$ for all ride segments and tries to reassign ride segments to drivers. For every driver route in the current solution, we check whether we can assign the ride segments to other drivers. A driver $k$ can be discarded if all ride segments in her route $\routeDriverCND_{k}(S)$ can be added to existing routes of the other drivers $l \in \bar{K}(S): l \neq k$.
This operator aims at improving the solution by eliminating drivers, which directly impacts the objective value. 

\paragraph{Postpone}
This operator aims to delay some of the departure times $t(S)$ of the current solution $S$, which might allow for more efficient driver schedules. Accordingly, this operator utilizes the time windows of our problem. After shifting the departure times of one or multiple rides, we recompute the driver routes $\routeDriverCND(S)$, while vehicle routes $\routeVehicleCND(S)$ remain unchanged. To update the driver routes, all operators use the procedure described in Section \ref{sss:constrHeuristic}.

\paragraph{Prepone}
Similar to the \textit{Postpone} operator, the \textit{Prepone} operator  aims at improving the current solution $S$ by shifting some of the departure times $t(S)$ and obtaining better driver routes $\routeDriverCND(S)$. While the \textit{Postpone} operator delays the departure times of selected rides, the \textit{Prepone} operator shifts the departure of selected rides to earlier points in time. Again, vehicle routes $\routeVehicleCND(S)$ remain unchanged.

\paragraph{Insert Intermediate Stop (Random)}
The previously presented operators focus on adjusting driver routes and departure times while fixing the vehicle routes. In contrast, the operator \textit{Insert Intermediate Stop (Random)} modifies the vehicle routes~$\routeVehicleCND(S)$ by inserting intermediate stops at service stations for potential driver exchanges. First, we choose a ride segment to insert an intermediate stop. Here, we prioritize segments with long durations as splitting these into shorter segments might allow for more flexible and thus more efficient driver routing. However, for the sake of diversification, we do not strictly select the ride segment with the longest travel duration but rather randomize the selection to a certain degree controlled by a perturbation factor~$\paramDeterminism$. Given a list of $n$ ride segments  ordered by their descending travel duration, we draw a random number $y \in [0,1)$ and select the ride segment at position $\lfloor y^{\paramDeterminism} n \rfloor$. Second, we randomly choose a service station reachable within the detour limit to insert into this ride segment, which changes some of the departure times~$t(S)$. We adjust the vehicle route accordingly. Finally, we update the driver routes~$\routeDriverCND(S)$. 

\paragraph{Insert Intermediate Stop (Shortest Detour)}
The operator \textit{Insert Intermediate Stop (Shortest Detour)} modifies vehicle routes by inserting intermediate stops at service stations. However, instead of selecting a service station randomly, we now choose a service station based on the additional driving time that the change implies. As in the previously described operator, we first select a ride segment. Here, we prioritize ride segments with long travel times. 
Next, we rank all service stations according to the descending detour they cause. From this list we do not simply select the station leading to the shortest detour but again randomize the selection as we do for the selection of ride segments.
Then, we adjust the vehicle routes $\routeVehicleCND(S)$ and departure times~$t(S)$ accordingly and finally update the driver routes~$\routeDriverCND(S)$.

\paragraph{Insert Intermediate Stop (Highest Synchronization Potential)}
As the two preceeding operators, this operator modifies vehicle routes by inserting intermediate stops at service stations. Here, we select service stations based on their synchronization potential, which we define through the number of nodes in routes of the current solution that share the same physical location. Given a list of service stations ordered by their descending synchronization potential, we again randomize the selection of service stations. Based on the selected service station, we adjust the vehicle routes $\routeVehicleCND(S)$ and departure times~$t(S)$ and update the driver routes~$\routeDriverCND(S)$.

\paragraph{Remove Intermediate Stop}
This operator aims at removing redundant stops at service stations to create more efficient driver routes and schedules. It modifies the vehicle routes $\routeVehicleCND(S)$ in a solution $S$ by eliminating intermediate stops until no further improvement is possible. Then, we recompute departure times~$t(S)$ and driver routes $\routeDriverCND(S)$ accordingly.

\section{Experimental Setup}\label{s:compStudy}
This section defines the experimental design for our computational study. In Section \ref{ss:instances}, we describe our real-world instances before we specify the experimental setting for our computational study in Section \ref{ss:paramFitting}.

\subsection{Instances}\label{ss:instances}
We derive our instances from a large real-world data set provided by Flixbus, a leading coach company in Europe.
This data set comprises more than 3,000 rides to be executed all over Europe on a single day in the summer high season of 2019. It includes all stops to be visited on each ride, scheduled departure times at bus stops, the travel duration for every ride segment, and the subcontractor responsible for executing the ride. Additionally, Flixbus provided us with the geocoordinates of service stations available for driver exchanges. To approximate the travel duration to and from service stations, we first calculate the orthodromic distance with the Haversine formula and then multiply this distance with a conversion factor provided by Flixbus.

We create instances from this data set as follows:
first, we transform the scheduled departure time $\tau_v$ of every customer node $v \in \setOfNodesFlat \backslash \{0\}$ into time windows $[\paramEarliestDep_{v}, \paramLatestDep_{v}]$ of width $\paramTimeWindowSize$, such that $\paramEarliestDep_{v} = \tau_v - \frac{\paramTimeWindowSize}{2}$ and
$\paramLatestDep_{v} = \tau_v + \frac{\paramTimeWindowSize}{2}$.
Based on the time window of a node and the level of time discretization $\paramIntervalLength$, we  define the set of discrete departure times as described in Section \ref{ss:graphRepr}. 
Next, we exclude service stations that cannot be visited within the detour limit $\paramDetourLimit$  for each ride.
Finally, as drivers may steer only buses of their company, we separate the data by subcontractor, which leads to 351 instances. These instances vary considerably in terms of size and complexity.
Therefore, we cluster them into three groups depending on the number of time-expanded arcs $\setOfArcsTimeExp$: i) small instances with $|\setOfArcsTimeExp| < 1,000$, ii) medium instances with $1,000 \leq |\setOfArcsTimeExp| < 5,000$, and iii) large instances with $5,000 \leq |\setOfArcsTimeExp|$.
Table \ref{tab:instanceAttrTd10} gives an overview of aggregated instance attributes for each group, displaying the number of rides ($|\setOfRides|$), the number of ride segments ($|\mathcal{RS}|$), the number of non-time-expanded nodes ($|\setOfNodesFlat|$), the number of customer locations ($\sum_{r \in \setOfRides} |\mathcal{C}_{r}|$), and the number of available service stations ($\sum_{r \in \setOfRides} | \setOfServiceStationsPerR |$).
For all instances, we consider a time window width of $\paramTimeWindowSize = 10$ minutes, a detour limit of $\paramDetourLimit = 10$ minutes, and a time discretization of $\paramIntervalLength = 10$ minutes, if not stated otherwise.
Note that we further decompose some of the instances into independent sub-instances. If a subcontractor operates multiple ride sets that do not share any bus stops or service stations, there is no possibility to transfer drivers between rides of different sets. Consequently, we optimize driver routes and schedules for these sets independently. 

\begin{figure}[tb]
	\begin{floatrow}
		\ttabbox{%
			\small
			\begin{tabular}{lrlrrrrrrr} 					
				\toprule							
				Instance size & \# Instances & & &
				$|\setOfRides|$ & $|\mathcal{RS}|$ &&
				$|\setOfNodesFlat|$ & $\sum_{r \in \setOfRides} |\mathcal{C}_{r}|$ & $\sum_{r \in \setOfRides} | \setOfServiceStationsPerR |$\\
				\midrule
				\multirow{4}{*}{Small} & \multirow{4}{*}{101} & Average && 2 & 13  && 73 & 15 & 58\\ 
				& & Median && 2 & 12  && 67 & 14 & 50\\
				& & Minimum && 1 & 4 && 13 & 5 & 2\\
				& & Maximum && 6 & 34  && 155 & 38 & 136\\
				\midrule
				\multirow{4}{*}{Medium} & \multirow{4}{*}{158} & Average && 6 & 34  && 274 & 40 & 233\\
				& & Median && 6 & 32  && 257 & 37 & 216\\
				& & Minimum && 2 & 11  && 55 & 14 & 10\\
				& & Maximum && 22 & 85  && 621 & 96 & 556\\
				\midrule
				\multirow{4}{*}{Large} & \multirow{4}{*}{92} & Average && 24 & 154  && 3,543 & 178 & 3,368\\
				& & Median && 17 & 111  && 1,044 & 127 & 867\\
				& & Minimum && 4 & 22  && 320 & 26 & 238\\
				& & Maximum && 168 & 1,153 && 74,130 & 1,321 & 74,013\\
				\bottomrule	
			\end{tabular}				
		}{%
\caption{Aggregated instance attributes}%
\label{tab:instanceAttrTd10}%
}
\end{floatrow}
\end{figure}

\subsection{Algorithm Design and Parameter Fitting}\label{ss:paramFitting}
We executed all numerical experiments in a single core setting on a cluster of 812 machines, each having 28 Intel Xeon E5-2690 v3 processors (2.6 GHz) with 32 GB RAM each, running on Linux. Both the \gls{MIP} and our matheuristic were implemented as a single core thread using the Gurobi Python Interface with Gurobi~9.1.1 and Python~3.8.8. In our matheuristic, we use Gurobi callbacks to alternate between the \gls{MIP} and the local search phase.
We set the time limit to $\paramTimeLimit$ seconds for all methods. 

We fitted the algorithmic parameters of our matheuristic based on a subset of 35 instances, including 10\% of all instances per instance size. Initial value sets for the algorithmic parameters have been identified during preliminary computations in the development phase.
In a systematic study, we tested these predefined sets by varying a single parameter at a time and fixing it to the value providing the best results. Table \ref{tab:paramSetting} 
\begin{figure}[]
	\begin{floatrow}
		\ttabbox{%
			\small
			\begin{tabular}{llcccclcccclccc} 					
				\toprule
				Parameter && \multicolumn{3}{c}{$\paramTimeLimitLBImpr$} &&& \multicolumn{3}{c}{$\paramTimeLimitLocalSearchCB$} &&& \multicolumn{3}{c}{$\paramSolutionsKept$}\\							
				\cmidrule{3-5}
				\cmidrule{8-10}
				\cmidrule{13-15}
				Values & 
				& 5\% & 10\% & \textbf{15\%} && 
				& 5\% & \textbf{10\%} & 15\% && 
				& \textbf{1\%} & 10\% & 20\% \\
				\cmidrule{3-15}
				$\Delta \objectiveValue$ [\%] & 
				&0.48 & 0.32 & 0.06 &&
				&0.45 & 0.06 & 0.17 &&
				&0.06 & 0.28 & 0.44 \\
				\midrule
				\midrule
				Parameter && \multicolumn{3}{c}{$\paramTimeLimitMIPsearch$} &&& \multicolumn{3}{c}{$\paramDeterminism$} &&& \multicolumn{3}{c}{$\paramSolutionsKeptAbsMin$}\\
				\cmidrule{3-5}
				\cmidrule{8-10}
				\cmidrule{13-15}
				Values & 
				& 20\% & \textbf{50\%} & 80\% && 
				& 3 & \textbf{5} & 10&& 
				& \textbf{10} & 15 & 20 \\
				\cmidrule{3-15}
				$\Delta \objectiveValue$ [\%] &
				&0.28 & 0.06 & 0.46 &&
				&0.46 & 0.06 & 0.46 &&
				&0.06 & 0.28 & 0.44 \\
				\bottomrule	
			\end{tabular}				
		}{%
			\caption{Parameter setting}%
			\label{tab:paramSetting}%
		}
	\end{floatrow}
\end{figure}

gives an overview of the parameter values considered and highlights final values in bold for
the time limit of the destructive bound improvement phase~($\paramTimeLimitLBImpr$), the time limit of the \gls{MIP} solver before injecting a start solution~($\paramTimeLimitMIPsearch$), the time limit of the local search phase during callbacks~($\paramTimeLimitLocalSearchCB$), the perturbation factor~($\paramDeterminism$), the share of solutions kept ($\paramSolutionsKept$), and the absolute minimum number of solutions kept~($\paramSolutionsKeptAbsMin$). The table further displays the average deviation over 10 runs from the best solution found with respect to the objective function value~($\Delta \objectiveValue$). 
Note that the best parameter setting results in an average deviation of 0.06\%, indicating that no global best setting for all instances exists.

\begin{figure}[bt]
	\begin{floatrow}
		\ffigbox[.8\textwidth]{%
			\includegraphics[width=.75\textwidth]{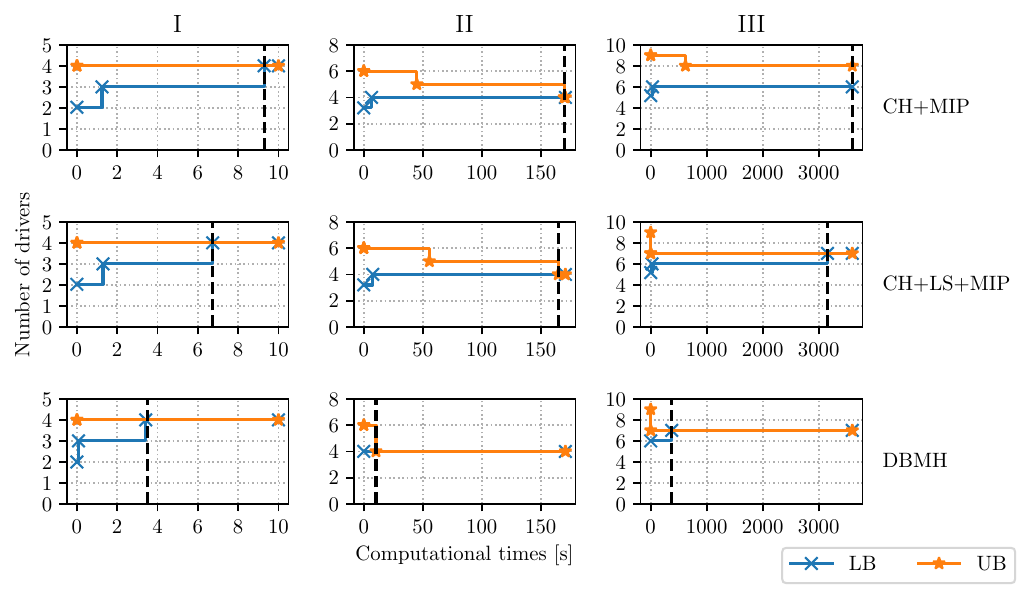}%
		}{%
			\caption{Convergence of CH+MIP, CH+LS+MIP, and \gls{DBMH}}%
			\label{fig:convergence}%
		}
	\end{floatrow}
\end{figure}

Next, we use the same subset of instances for our algorithm design. In particular, we are interested in the performance gain resulting from enhancing a basic matheuristic, such as a warm-started \gls{MIP}, with a \gls{LS} component providing an improved start and the \gls{DBI} component lifting the lower bound.
Therefore, we compute 10~runs for every instance and compare the results of three different matheuristic variants: i)~the \gls{MIP} warm-started with the \gls{CH} solution (\gls{CH}+\gls{MIP}), ii)~~the \gls{MIP} warm-started with the \gls{CH}+\gls{LS} solution (\gls{CH}+\gls{LS}+\gls{MIP}), and iii)~the full \gls{DBMH} as described in Section \ref{ss:solApproach}. 
Over all instances, warm-starting the \gls{MIP} with the \gls{CH}+\gls{LS} start solution reduces the number of drivers by up to 7 compared to the \gls{CH} solution and increases the number of instances solved to optimality within the time limit by 6.54\%. Moreover, the computational time for instances solved to optimality can be reduced by 23.63\% or 51.33~seconds on average. 
The results further show that the lower bound increases much faster when integrating the \gls{DBI}, leading to better convergence and reducing computational times for instances solved to optimality by 54.87\% or 169.44~seconds on average. 
We illustrate the impact of the \gls{LS} and the \gls{DBI} components in Figure~\ref{fig:convergence},

which plots the development of lower and upper bounds for three exemplary instances of small and medium size. 
The plots for instance~III show that an improved start solution is particularly effective with a tight lower bound. While it takes the \gls{CH}+\gls{LS}+\gls{MIP} method more than 3,000~seconds to lift the lower bound and prove optimality of the \gls{CH}+\gls{LS} start solution, the \gls{DBMH} needs less than 400~seconds to terminate. 
Furthermore, the \gls{DBI} cannot only prove optimality of heuristic start solutions but also finds optimal solutions while checking for feasibility of a particular objective function value, as the plots for instance~II demonstrate. Here, the \gls{DBI} starts from a constructive lower bound~$cLB=4$ and  finds a feasible and hence optimal solution with an objective function value of 4 while checking for feasibility of $LB=UB=4$. Over all test runs, the \gls{DBMH} solves 5.26\% more instances to optimality than the \gls{CH}+\gls{LS}+\gls{MIP} method.
As expected, tight lower and upper bounds improve convergence, significantly reducing computational times for instances solved to optimality by 70.14\%. Therefore, we chose the \gls{DBMH} with advanced start and destructive bound improvement over the basic warm-started \gls{MIP} variants. 

\section{Results}\label{s:results}
This section presents the results of our computational study. 
We first evaluate the performance of our approach in Section~\ref{ss:MIPvsHybridvsCND}. In Section~\ref{ss:algPerformance}, we analyze the performance of the individual algorithmic components. Finally, we derive managerial insights in Section~\ref{ss:manInsights}.

\subsection{Computational Performance}\label{ss:MIPvsHybridvsCND}

To validate the performance of our \gls{DBMH}, we compare its results to solutions of the \gls{MIP} and to heuristic solutions, considering all 351 real-world instances. We generate heuristic solutions by applying the \gls{CH} first and the \gls{LS} afterward. 
We conducted 10 runs per instance for every solution approach.  
Table~\ref{tab:resultsMipHybridCNDtd10} provides an overview of the aggregated results by instance size and solution method. 
%
\begin{figure}[]
	\begin{floatrow}
		\ttabbox{%
			\small
			\begin{tabular}{lllrrrrr} 					
				\toprule
				\multicolumn{2}{l}{Instance size}&& Small & Medium & Large && All\\
				\multicolumn{2}{l}{\# Instances}&& 101 & 158 & 92 && 351\\
				\multicolumn{2}{l}{Average best $\objectiveValue$}&& 2.63 & 6.58 & 29.75 && 11.52\\
				\midrule
				\multirow{4}{*}{MIP} &  opt.solved[\%] & &99.55 & 51.25 & 1.60 && 53.02\\
				& t[s] & & 26.43 & 312.98 & 1,849.17 && 162.34 \\
				& gap[\%]  & & 0.00 & - & - & & - \\ 
				& $\Delta \objectiveValue$[\%] & & 0.00 & - & - & & - \\
				& $\sum \objectiveValue$ & & 266 & - & - & & - \\
				\midrule
				\multirow{4}{*}{\gls{DBMH}} &  opt.solved[\%] & & 100.00 & 60.00 & 4.26 && 57.69 \\
				& t[s] 					& & 5.44 & 422.13 & 709.00 && 204.53\\
				& gap[\%]  				& & 0.00 & 6.81 & 32.12 & & 11.77 \\ 
				& $\Delta \objectiveValue$[\%] & & 0.00 & 0.24 & 0.06& & 0.12\\
				& $\sum \objectiveValue$ & & 266 & 1,042 & 2,739 & & 4,047 \\
				\midrule
				\multirow{4}{*}{CH+LS} &  opt.solved[\%] & & 68.70 & 24.43 & 2.15 && 32.09\\
				& t[s] & & 0.02 & 0.45 & 97.57 &&24.39\\
				& gap[\%]  & & 8.33 & 17.45  & 33.58  & & 19.23\\
				& $\Delta \objectiveValue$[\%] & & 11.61 & 15.10 & 1.60 & & 10.39 \\
				& $\sum \objectiveValue$ & & 300 & 1,157 & 2,745 & & 4,202 \\
				\bottomrule	
			\end{tabular}				
		}{%
			\caption{Aggregated results on real-world instances}%
			\label{tab:resultsMipHybridCNDtd10}%
		}
	\end{floatrow}
\end{figure}

The table header states the number of instances (\# Instances) and the average best driver count per instance (Average best $\objectiveValue$) for each instance size. In the table's main body, we report the percentage of instances solved to optimality (opt.solved[\%]), the average computational time to optimality (t[s]), the average gap between the best solution and the best bound (gap[\%]),  the average deviation from the best-known solution ($\Delta \objectiveValue$[\%]), and the accumulated driver count over all instances per instance size ($\sum \objectiveValue$) for every solution method. 
The results in Table~\ref{tab:resultsMipHybridCNDtd10} demonstrate that our \glsfirst{DBMH} outperforms the \gls{MIP} and the \gls{CH}+\gls{LS} on real-world instances.
Specifically, on medium and large instances, our matheuristic solves more instances to proven optimality than the \gls{MIP}. This result demonstrates the relevance of the heuristic components and the destructive bound improvement for obtaining optimal solutions. 
Moreover, these additional components allow us to reduce computational times by 79.42\% for small instances compared to the \gls{MIP}.
Another advantage of our approach is its ability to provide feasible solutions for all instances within seconds. By contrast, the \gls{MIP} fails to find any feasible solution for the majority of medium and large instances, which is why we cannot state the \gls{MIP}'s average gap, deviation from best-known solutions, and total driver count for these instances.
While the combination of \gls{CH} and \gls{LS} performs best with respect to computational times, its solutions deviate by 10.39\% on average from best-known solutions and require a total of $155$ additional drivers compared to the \gls{DBMH}, indicating the importance of the destructive bound improvement and the \gls{MIP} to reach optimal solutions. 
For large instances, however, the deviation of heuristic solutions from best-known solutions is relatively low with 1.60\% as our matheuristic often terminates without finding the optimal solution in these cases.

Recalling the requirements resulting from the practical application described at the beginning of Section~\ref{ss:solApproach}, the results in Table~\ref{tab:resultsMipHybridCNDtd10} indicate that our \glsfirst{DBMH} is the most appropriate solution method for practical applications as it generates feasible solutions quickly and provides information about the obtained solution quality.

\subsection{Performance of the Algorithmic Components}\label{ss:algPerformance}
We now analyze the impact of the components of our \glsfirst{DBMH} on solution quality and computational time in detail. 
%
\begin{figure}[t]
	\begin{floatrow}
		\ttabbox[]{%
			\small
			\begin{tabular}{llrr} 					
				\toprule
				&& \multicolumn{2}{c}{$\Delta \; LB [\%]$}\\
				\cmidrule{3-4}
				Instance size&& avg & med\\
				\midrule
				Small	&&	26.55 & 0.00\\
				Medium	&&	17.58 & 16.67\\
				Large	&&	2.99 &  0.00\\
				\midrule
				All && 20.14 & 0.13\\
				\bottomrule	
			\end{tabular}
		}{%
			\caption{Impact of destructive bound improvement}%
			\label{tab:constructiveVsDestructiveLB}%
		}
	\end{floatrow}
\end{figure}

\CenterFloatBoxes
\begin{figure}[]
	\begin{floatrow}
		\ffigbox[.5\textwidth]{%
			\includegraphics[width=.45\textwidth]{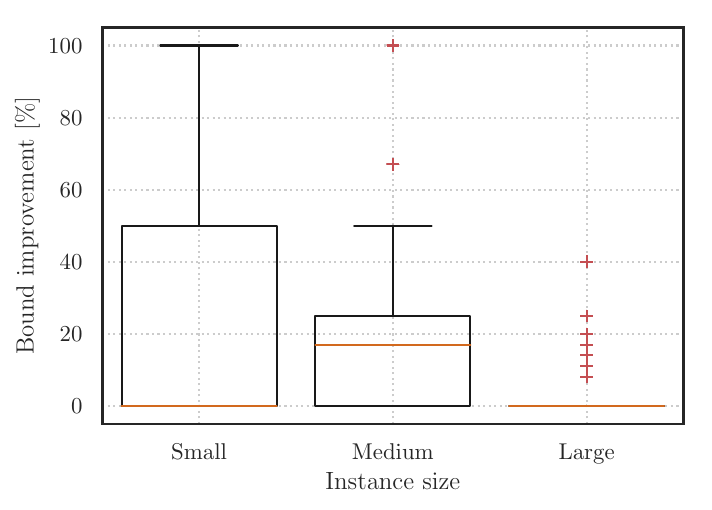}%
		}{%
			\caption{Relative bound improvement}%
			\label{fig:boundImpr}%
		}
		\ffigbox[.5\textwidth]{%
			\includegraphics[width=.45\textwidth]{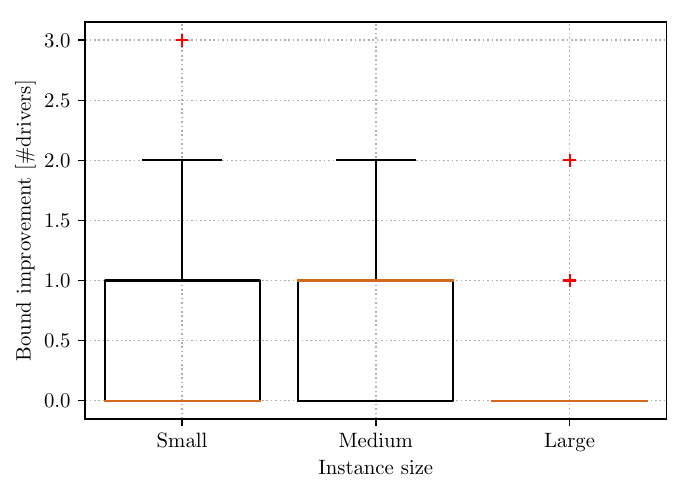}%
		}{%
			\caption{Absolute bound improvement}%
			\label{fig:boundImprAbs}%
		}
	\end{floatrow}
\end{figure}

\paragraph{Impact of Destructive Bound Improvement}
First, we study the impact of the \gls{DBI} on the lower bound considering all 351 real-world instances.
Table~\ref{tab:constructiveVsDestructiveLB} depicts the average (avg) and median (med) improvement of the \gls{dLB} compared to the \gls{cLB}, $\Delta \; LB \; = \frac{dLB-cLB}{cLB} \; [\%]$, grouped by instance size. 
Figures~\ref{fig:boundImpr} and \ref{fig:boundImprAbs} contain additional information on the distribution of the bound improvement.
By definition, $dLB \geq cLB$ holds. As can be seen, the \gls{DBI} increases the lower bound significantly and by 20.14\% on average. The \gls{dLB} is particularly effective for small and medium instances. The large delta between average and median improvement as well as the boxplots in Figures~\ref{fig:boundImpr} and \ref{fig:boundImprAbs} show that the level of improvement varies notably depending on the instance size. 
While the \gls{DBI} achieves significant bound improvements of up to 100\% or three drivers for some instances, there exist many small and large instances for which the \gls{DBI} cannot improve the \gls{cLB}. Recalling the results from Table~\ref{tab:resultsMipHybridCNDtd10}, there are two reasons for this observation: first, trivial instances can already be solved optimally through the \gls{CH}, the \gls{LS}, and the \gls{cLB}, and hence no further bound improvement is possible. Second, large instances may bear a complexity for which the \gls{DBI} is not able to improve the bound in the given time limit.

\begin{figure}[]
\begin{floatrow}
\ttabbox[.5\textwidth]{%
	\small
  \begin{tabular}{llrrrrr} 					
  	\toprule
  	&&\multicolumn{5}{c}{$\Delta \objectiveValue$[\%]}\\
  	\cmidrule{3-7}
  	Instance size && \gls{cLB} & \gls{dLB} && \gls{LS} & \gls{CH} \\
  	\midrule
  	Small	&&	-15.70  	& 0.00	&&	11.97 	& 15.27 \\
  	Medium	&&	-18.93  	& -3.47	&&	18.59 	& 24.40 \\
  	Large	&&	-17.78 	& -5.56	&&	31.11 	& 31.11 \\
  	\midrule
  	All && 	-17.23 	& -1.65 	&&	15.15 	& 19.62 \\
  	\bottomrule	
  \end{tabular}
}{%
  \caption{Average quality of bounds and start solutions}%
  \label{tab:boundTightness}%
}
\ffigbox[.5\textwidth]{%
	\includegraphics[width=.45\textwidth]{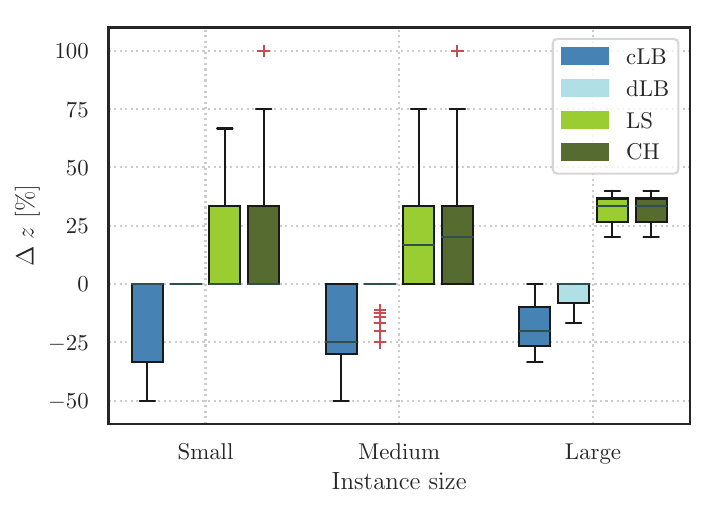}%
}{%
	\caption{Deviation from optimal objective value}%
	\label{fig:boundTightness}%
}
\end{floatrow}
\end{figure}

\paragraph{Quality of Bounds and Start Solutions}
Next, we analyze the tightness of both lower bounds as well as the quality of the constructive and the improved start solution. In this context, we analyze the average deviation over 10 runs from the optimal solution ($\Delta \objectiveValue [\%]$) for the \gls{cLB}, the \gls{dLB}, the initial solution obtained by the \gls{CH}, and the improved initial solution after applying the \gls{LS}. 
Table~\ref{tab:boundTightness} depicts aggregated results grouped by instance size, while Figure~\ref{fig:boundTightness} shows boxplots of the deviation for every instance size.
Note that for this and the remaining analyses in this section, we only use instances solved to optimality, which include 101 small, 95 medium, and 3 large instances.
\CenterFloatBoxes

Figure~\ref{fig:boundTightness} shows that the quality of the lower bounds and the start solutions varies notably among different instances, in particular for small instances. While the \gls{cLB}, the \gls{dLB}, the \gls{CH}, and the \gls{LS} match the optimal objective function value for more than half of the small instances, there exist some instances for which the \gls{cLB} deviates up to -50\% and the start solutions deviate up to 100\% from the optimal objective function value.
Only the \gls{dLB} is generally tight, deviating by -1.65\% on average from the optimal solution as the results in Table~\ref{tab:boundTightness} and Figure~\ref{fig:boundTightness} show. In absolute terms, the \gls{cLB} and \gls{dLB} deviate from the optimal driver count by up to $-3$ and $-2$ drivers, respectively.
The heuristic start solutions are far from optimal solutions, with a $\Delta \objectiveValue$ per instance of up to $6$ drivers for the \gls{CH} and up to $4$ drivers for the \gls{LS}, which is in line with the results in Table~\ref{tab:resultsMipHybridCNDtd10}. 
Overall, the average quality of the bounds and the start solutions decreases with growing instance size.

\paragraph{Generation and Proof of Optimal Solutions}
Next, we investigate which algorithmic components generate and prove optimal solutions. 
We differentiate four stages in which optimality can be found: 
i) after generating an initial solution by applying the \gls{CH} and the \gls{LS}\gls{init}, 
ii) after applying the \gls{LS} in the \gls{CB}, 
iii) during the \gls{DBI},  
or iv) while solving the \gls{MIP}.	
Table~\ref{tab:resultsSoluntionFound} shows the percentage share of instances solved to optimality in each stage. 
\begin{figure}[t]
	\begin{floatrow}
		\ttabbox[\textwidth]{%
			\small
			\begin{tabular}{lcccccc} 					
				\toprule							
				& & \multicolumn{2}{c}{Heuristic components} & & \multicolumn{2}{c}{MIP components} \\
				\cmidrule{3-4}
				\cmidrule{6-7}
				Instance size & & \gls{init} & \gls{CB} & & \gls{DBI} & \gls{MIP}\\
				\midrule
				Small & & 68.40\%  & - &  & 31.60\% & - \\
				Medium & & 40.16\%  & 0.84\% &  & 46.58\% & 12.42\% \\
				Large & & - & - && 66.67\% & 33.33\% \\
				\midrule
				All & & 54.88\% & 0.39\%  & & 38.77\% & 5.96\% \\
				\bottomrule	
			\end{tabular}
		}{%
			\caption{Components finding optimal solutions}%
			\label{tab:resultsSoluntionFound}%
		}
	\end{floatrow}
\end{figure}

The results in Table~\ref{tab:resultsSoluntionFound} demonstrate why our matheuristic outperforms the \gls{MIP} and combination of \gls{CH} and \gls{LS}. 
While the heuristic components can find optimal solutions in particular for small instances, \gls{MIP} components, primarily the \gls{DBI}, are essential to find and prove optimal solutions for larger instances. 
Specifically, the heuristic start solution is optimal for 54.88\% of all instances, which explains the low computational times of our solution approach for small instances. 
Although the heuristic components cannot provide optimal solutions for large instances, they are still important for providing feasible solutions. For the majority of instances not solved to optimality within the time limit, the heuristic start solution cannot be improved any further by the \gls{MIP}. This explains the rather low deviation of the \gls{CH}+\gls{LS} solutions from best-known solutions for large-size instances in Table~\ref{tab:resultsMipHybridCNDtd10}.

\paragraph{Distribution of Computational Time}
Finally, we analyze the distribution of computational time among the different algorithmic components. Due to large differences in absolute computational times, we show the percentage of the total computational time for every component and instance. Table~\ref{tab:runtimes_opt} shows the aggregated computational time shares grouped by component and the average computational time per instance in seconds (t[s]) per instance size.
\begin{figure}[]
	\begin{floatrow}
		\ttabbox[\textwidth]{%
			\small
			\begin{tabular}{lrrrrrrr} 					
				\toprule							
				& && \multicolumn{2}{c}{Heuristic components} & & \multicolumn{2}{c}{MIP components} \\
				\cmidrule{4-5}
				\cmidrule{7-8}
				\parbox{1cm}{Instance size} & t[s] && \gls{init} & \gls{CB} && \gls{DBI} & \gls{MIP} \\
				\midrule
				Small	& 7.20		&& 48.55\%  	& 0.00\% 	&& 51.45\%	& 0.00\% \\
				Medium	& 448.68	&&	16.31\%  	& 0.01\% 	&&	70.46\%	& 13.22\% \\
				Large	& 669.23	&&	0.27\% 		& 0.15\% 	&&	80.22\%	& 19.36\% \\
				\midrule
				All & 217.38	&& 	33.20\% 	& 0.01\% 	&&	60.50\% 	& 6.29\% \\
				\bottomrule	
			\end{tabular}
		}{%
			\caption{Computational time distribution}%
			\label{tab:runtimes_opt}%
		}
	\end{floatrow}
\end{figure}

As the results in Table~\ref{tab:runtimes_opt} show, the \gls{MIP} components account for 66.79\% of the total computational time on average. 
For small instances, the computational times are evenly distributed between creating the start solution and the \gls{DBI}.
With growing instance size, the percentage share of computational time spent in heuristic components decreases due to two reasons: first, the growth rate of total computational times of the heuristic components is lower than the growth rate of \gls{MIP} components. Second, as reported in Table~\ref{tab:resultsSoluntionFound}, heuristic components are less successful in providing optimal or near-optimal solutions for large instances, and hence more time is spent in the \gls{MIP} components.
Also for instances not solved to optimality, solving the \gls{MIP} accounts for the major part of the computational time. In contrast, the time for the \gls{DBI} converges to its time limit $\paramTimeLimitLBImpr$. Computational times of the heuristic components add up to less than 1\% on average for instances not solved to optimality.\\\\
To conclude, the combination of heuristic elements with mathematical programming techniques enriched by a destructive bound improvement procedure enables our matheuristic to quickly solve real-world instances to optimality and to provide feasible solutions even for large instances within seconds.

\subsection{Managerial Insights}\label{ss:manInsights}
To derive managerial insights, we first compare our results to solutions currently obtained in practice and investigate the advantages of jointly considering all lines of one subcontractor compared to an isolated planning per line. Second, we evaluate the benefit of driver exchanges and intermediate stops. Finally, we study the impact of time windows, detour limits, and time discretization on solution quality and computational effort.

\begin{figure}[]
	\begin{floatrow}
		\ttabbox[\textwidth]{%
			\small
			\begin{tabular}{llrrrrrrr} 					
				\toprule
				\multirow{2}{1cm}{Instance size} && \multicolumn{3}{c}{Line} && \multicolumn{3}{c}{Subcontractor}\\
				\cmidrule{3-5}
				\cmidrule{7-9}
				&& avg & med & $\Delta \sum \objectiveValue$ && avg & med & $\Delta \sum \objectiveValue$\\
				\midrule
				Small &&8.87\%& 0.00\% & 35 && 10.28\% & 0.00\% & 41  \\
				Medium &&13.50\%& 12.50\% & 73 && 15.85\% & 16.67\% & 85\\
				Large 	&& 22.62\%& 25.00\% & 5 && 27.38\% & 28.57\% & 6 \\
				\midrule
				All&& 11.18\% & 0.00\% & 113 && 13.07\% & 5.56\% & 132\\
				\bottomrule	
			\end{tabular}
		}{%
			\caption{Savings compared to the status quo used in practice}%
			\label{tab:resultsCompLine}%
		}
	\end{floatrow}
	\vspace{.2cm}
	\begin{floatrow}
		\ffigbox[\textwidth]{%
			\includegraphics[width=.45\textwidth]{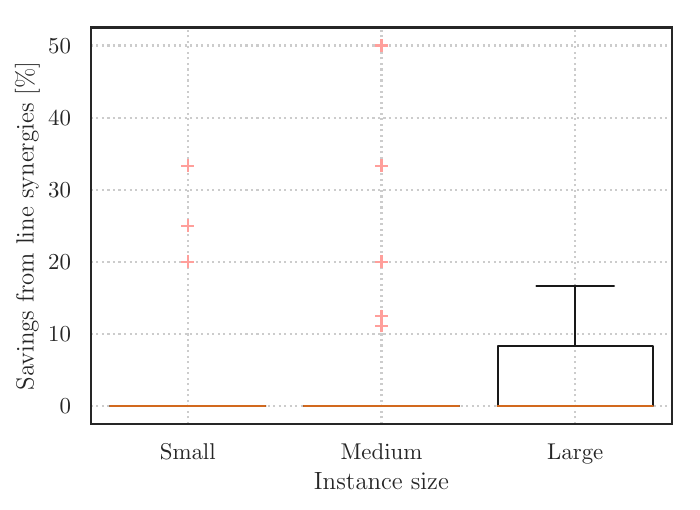}%
		}{%
			\caption{Savings from synergies between lines}%
			\label{fig:synergies}%
		}
	\end{floatrow}
\end{figure}

\paragraph{Status Quo Comparison}
In practice, driver assignments are planned manually by a team of network planners, each being responsible for a set of lines. Due to the size and complexity of instances on subcontractor level and the lack of sufficient software tools, drivers are scheduled for every line separately. This means that the network planner applies a similar logic as used in the construction heuristic described in Section~\ref{sss:constrHeuristic}. Therefore, we divide all 351 instances into subinstances at the line level and apply the construction heuristic to approximate manual solutions. 
To compare our methodology against this status quo benchmark, we apply our \glsfirst{DBMH} first to line-based instances and second to subcontractor-based instances. We then compare both results to the approximated manual results. Table~\ref{tab:resultsCompLine} shows the average~(avg) and median~(med) savings per instance as well as the total driver count reduction~($\Delta \sum \objectiveValue$) for both comparisons grouped by instance size and instance type.
The results in Table~\ref{tab:resultsCompLine} show that our matheuristic outperforms the status quo benchmark. Over all line-based instances, we can improve manual solutions by 11.18\% per instance on average, obtaining larger savings of 22.62\% with growing instance size. 
When considering subcontractor-based instances, we can further increase the average savings to 13.07\%. 
This result indicates that the primary savings result from the algorithm itself and only a minor part results from existing synergies between lines of the same subcontractor. 
Note that, in particular for small and medium instances, most instances do not have relevant synchronization potentials between lines, see Figure~\ref{fig:synergies}.
\CenterFloatBoxes

Only a small number of instances shows synergies between lines and hence considerable savings of up to 50\%, adding up to a reduction of 19 drivers over all instances.

In summary, our \glsfirst{DBMH} is able to improve the results currently achieved in practice. To reach globally optimal solutions, coach companies should optimize driver routes and schedules on subcontractor level. However, considering subcontractor-based instances is beneficial only for some instances and increases the problem complexity significantly.

\paragraph{Benefit of Driver Exchanges and Intermediate Stops}
A unique characteristic of the \gls{DRSPMS} is the consideration of intermediate stops for driver exchanges. On ride segments with a duration of more than $\paramContSteering$~hours, a stop at a service station is required to make this segment feasible. However, on all other ride segments, stops at service stations are optional but may enable more efficient operations. 
While a large number of service stations for driver exchanges might increase the potential for cost savings, it also increases the problem complexity. 
In this context, we analyze the savings potential from driver exchanges in general and from considering driver exchanges at intermediate stops in addition to driver exchanges at regular stops. Here, we consider only instances solved to optimality. To get a sufficiently large set of optimally solved instances for these analyses, we increase the time limit for medium and large instances to $2$ and $4$ hours, respectively. Furthermore, we study the impact of intermediate stops on computational efficiency.

First, we compare solutions with driver exchanges at intermediate stops to solutions without driver exchanges. Here, we consider the possibility of double driver teams for rides with a duration exceeding the continuous steering time limit $\paramContSteering$. It is worth mentioning that without the possibility of driver exchanges, 176 of the 351 instances are infeasible as they include rides with a travel duration higher than the daily driving time of a double-driver team. 
Table~\ref{tab:impactDE} shows the average and median savings per instance as well as the total driver count reduction~($\Delta \sum \objectiveValue$) from considering driver exchanges grouped by instance size. Note that we only take instances into account that were solved optimally in both settings, such that our comparison bases on 149 instances in total
, including 82 small, 60 medium, and 7 large instances.
The results reported in Table~\ref{tab:impactDE} demonstrate the importance of driver exchanges for obtaining cost-efficient driver routes and schedules. 
Introducing the possibility of driver exchanges at regular and intermediate stops improves the solution quality by 42.69\% on average compared to solutions without driver exchanges. Over all 149 instances, the consideration of driver exchanges reduces the total number of required drivers by 406. Note that also for instances not solved to optimality within the time limit, the solutions with driver exchanges consistently improve solutions without driver exchanges.

%
\begin{figure}[]
	\begin{floatrow}
		\ttabbox{%
			\small
			\begin{tabular}{llrrr} 					
				\toprule
				Instance size && Average & Median & $\Delta \sum \objectiveValue$\\
				\midrule
				Small && 39.87\% &  50.00\% & 124\\
				Medium && 46.88\% &  50.00\% & 244 \\
				Large && 39.73\% & 50.00\% & 38 \\
				\midrule
				All && 42.69\% &  50.00\% & 406 \\
				\bottomrule	
			\end{tabular}
		}{%
			\caption{Savings from driver exchanges}%
			\label{tab:impactDE}%
		}
	\end{floatrow}
\end{figure}

\CenterFloatBoxes
\begin{figure}[]
	\begin{floatrow}
		\ttabbox[.45\textwidth]{%
			\small
			\begin{tabular}{llrrr} 					
				\toprule
				&&  \multicolumn{3}{c}{\gls{DE-RIS} vs. \gls{DE-RS}}\\
				\cmidrule{3-5}
				Instance size && Average & Median & $\Delta \sum \objectiveValue$\\
				\midrule
				Small && 1.13\% & 0.00\% & 5\\
				Medium && 2.13\% & 0.00\% & 14\\
				Large && 1.11\% & 0.00\% & 1\\
				\midrule
				All && 1.64\% & 0.00\% & 20\\
				\bottomrule	
			\end{tabular}
		}{%
			\caption{Savings from intermediate stops}%
			\label{tab:impactDE-IS}%
		}
		\hspace{.05\textwidth}
		\ffigbox[.45\textwidth]{%
			\includegraphics[width=.45\textwidth]{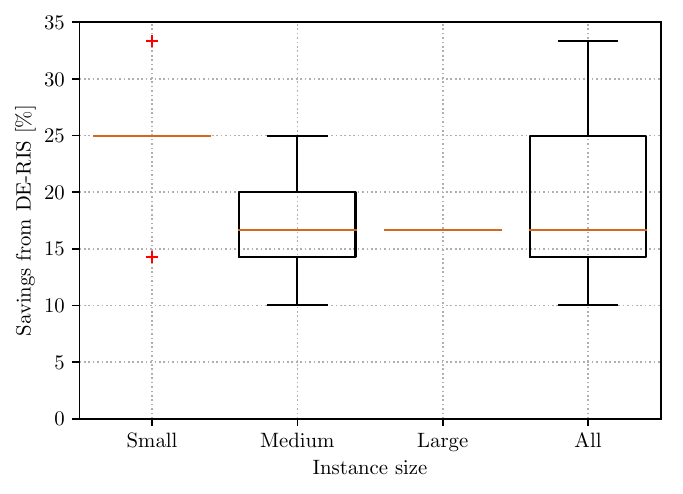}%
		}{%
			\caption{Distribution of savings from intermediate stops}%
			\label{fig:benefitIS}%
		}
	\end{floatrow}
\end{figure}

Next, we evaluate the savings potential from considering driver exchanges at intermediate stops in addition to driver exchanges at regular stops. Therefore, we compare optimal driver counts for the problem with \glsfirst{DE-RIS}~(\gls{DE-RIS}) to optimal driver counts for the problem with \glsfirst{DE-RS}~(\gls{DE-RS}). To compute solutions for the problem without intermediate stops, we set the detour limit~$\paramDetourLimit$ to zero and apply our \gls{DBMH}. Note that 14 out of 351~instances are infeasible in the \gls{DE-RS} setting due to ride segments exceeding the continuous steering time limit $\paramContSteering$. Only when allowing driver exchanges at intermediate stops all 351~instances are feasible. Table~\ref{tab:impactDE-IS} shows the average and median savings per instance as well as the total driver count reduction~($\Delta \sum \objectiveValue$) grouped by instance size. Our comparison includes all 220~instances solved to optimality in both settings, that is 100 small, 112 medium, and 8 large instances.
For 19~instances (5 small, 13 medium, 1 large), the consideration of intermediate stops results in improved solutions and reduces the total driver count by 20~drivers or 17.54\%. Figure~\ref{fig:benefitIS} plots the distribution of savings for these instances. For the remaining instances, the driver count difference is zero, which explains the comparatively low average savings per instance of 1.64\% reported in Table~\ref{tab:impactDE-IS}. Over all 220~instances solved to optimality, the option of driver exchanges at intermediate stops reduces the total driver count from~897 to 877~drivers, corresponding to a decrease of~2.23\%. 

Finally, we analyze the impact of intermediate stops on problem complexity and computational efficiency. Here, we consider all 337 instances feasible with and without intermediate stops. The consideration of intermediate stops offers additional flexibility and the potential to reduce the number of drivers required but also substantially increases the problem complexity and, thus, the computational effort. Specifically, computational times to optimality are 745.67 seconds higher on average than when considering driver exchanges only at regular bus stops. Furthermore, without intermediate stops, 55 more instances can be solved to proven optimality within the given time limit. Moreover, as the problem with intermediate stops is harder to solve, we might obtain worse solutions for some instances within the given time limit than we do for the problem without intermediate stops within the same time limit. This is the case for 9.20\% of the instances. 
Recall, however, that solutions for the problem without intermediate stops are also feasible for the problem with intermediate stops, which implies that we can always use the solution without intermediate stops for the problem with intermediate stops if the former is better.
Nevertheless, developing a method that can effectively and consistently exploit the potential of intermediate stops for large instances is an important step for future research. In this regard, adapting existing column generation approaches for integrated crew scheduling problems with retiming options, such as the methodology proposed by \cite{KliewerAmbergEtAl2012} or \cite{BachDollevoetEtAl2016} to include intermediate stops, seems promising.

In conclusion, driver exchanges at regular customer stops are essential to obtain feasible and cost-efficient driver routes and schedules. For some instances, the additional consideration of intermediate stops leads to further cost savings. Moreover, intermediate stops are required to ensure feasibility in some cases.

\paragraph{Impact of Different Detour Limits and Time Windows}
In general, the number of intermediate stops considered and thus the computational effort depend on two main factors: the detour limit $\paramDetourLimit$ and the departure time window size $\paramTimeWindowSize$, which are closely related as described in Section~\ref{ss:modelingConcept}.
Consequently, we jointly examine the effect of different detour limits and time window sizes. Table~\ref{tab:impactTWdetour} shows the average deviation from best-known solutions for different detour limits ($\paramDetourLimit [min]$) and time window sizes ($\paramTimeWindowSize [min]$). To evaluate the impact on the solution quality, we only consider 135 instances solved to optimality in all settings. Furthermore, we do not differentiate by instance size as the size depends on the number of arcs, which differ with time window width and detour limit.
The results in Table~\ref{tab:impactTWdetour} demonstrate that the configuration providing the largest degree of freedom (i.e., $\paramTimeWindowSize = 30$ minutes and $\paramDetourLimit = 30$ minutes) yields the best results. This degree of freedom reduces the total driver count by 6 compared to a rather restricted setting with $\paramTimeWindowSize = 10$~minutes and $\paramDetourLimit = 5$~minutes. 
Note that the impact of the time window size~$\paramTimeWindowSize$ is higher than the impact of the detour limit~$\paramDetourLimit$. Specifically, if the time window is sufficiently large, the detour limit can be set to a small value without any loss in solution quality. A higher detour limit can only achieve better results for narrow time windows. However, higher detour limits and wider time windows come at the expense of higher computational times. Increasing the time window size~$\paramTimeWindowSize$ from 10~minutes to 30~minutes increases the computational effort by more than a factor 100 on average. The impact of the detour limit on computational times is lower. A change of the detour limit~$\paramDetourLimit$ from 5~minutes to 30~minutes yields a computational time increase of less than factor~10.
Consequently, a setting that balances solution quality and computational times (e.g., $\paramTimeWindowSize = 10$~minutes and $\paramDetourLimit = 5$~minutes) should be chosen for practical applications.
Note that large time windows are generally unfavorable in practice because, in addition to driver routing and scheduling, many other planning problems, e.g., bus assignments, depend on the target departure times.

\paragraph{Impact of Different Time Discretization Levels}

Finally, we examine the impact of time discretization on solution quality and computational effort. The granularity of the time discretization in our underlying network determines the solution quality.
Increasing the level of detail and thus the degree of freedom in departure times increases the computational effort but might allow for more efficient driver routes. 
To study this effect, we compute results for all 351 instances considering a time interval length of $\paramIntervalLength = 5$ minutes and compare them to our results considering a time interval length of $\paramIntervalLength = 10$ minutes. 
Table~\ref{tab:impactTD} gives an overview of 
\begin{figure}
	\begin{floatrow}
		\ttabbox[.6\textwidth]{%
			\small
			\begin{tabular}{lrrrrrr} 					
				\toprule
				&& \multicolumn{3}{c}{$\paramTimeWindowSize \; [min]$}&&\\
				\cmidrule{3-5}
				&& $10$ & $20$ & $30$&&Avg\\
				\midrule
				\multirow{4}{1.5cm}{$\paramDetourLimit \; [min]$} & 5 & 1.60\% & 0.25\% & 0.00\% && 0.62\% \\
				&10 & 1.50\% & 0.25\% & 0.00\%&& 0.60\% \\
				&20 & - & 0.25\% & 0.00\% && 0.12\% \\
				&30 & - & - & 0.00\% && 0.00\%\\
				\midrule
				&Avg & 1.55\% & 0.25 \%& 0.00\%&&  -\\
				\bottomrule	
			\end{tabular}	
		}{%
			\caption{Impact of different detour limits and time windows}%
			\label{tab:impactTWdetour}%
		}
	
		\ttabbox[.4\textwidth]{%
			\small
			\begin{tabular}{llrr} 					
				\toprule
				&&\multicolumn{2}{c}{$\paramIntervalLength$ [min]}\\
				\cmidrule{3-4}
				&& 5 & 10 \\
				\midrule
				opt.solved [\%] && 50.86 & \textbf{57.69}\\
				$\Delta \objectiveValue^*$[\%] && \textbf{0.00} & 0.48\\
				$\Delta \objectiveValue$ [\%] && 2.11 & \textbf{0.32}\\
				t[s] && 199.13 & \textbf{18.53}\\
				\bottomrule	
			\end{tabular}
		}{%
			\caption{Impact of time discretization}%
			\label{tab:impactTD}%
		}
	\end{floatrow}
\end{figure}

the share of instances solved to optimality (opt.solved[\%]), 
the average deviation from best-known solutions for instances solved optimally in the respective setting ($\Delta \objectiveValue^*$[\%]), the average deviation from best-known solutions for all instances ($\Delta \objectiveValue$[\%]), and the average computational time ($t[s]$) per instance. 
As can be seen, the increased problem complexity outweighs the potential for better solution quality. 
While this increased level of detail leads to improved optimal solutions for three instances with a reduction of one driver each, the solution quality for instances not solved to optimality is worse with a time discretization of $\paramIntervalLength=5$ so that the total driver count increases by 45. Furthermore, the computational time to optimality is higher by a factor of 11 on the same instances, and the share of instances solved to optimality within the time limit decreases by 11.84\%.
Consequently, for practical applications that aim for a fair balance between solution quality and computational effort, we consider a time discretization of $\paramIntervalLength = 10$ minutes to be sufficient.

\section{Conclusion and Outlook}\label{s:conclusion}
In this paper, we presented a new type of driver routing and scheduling problem motivated by a practical application in long-distance bus networks. A unique characteristic of this problem is the possibility of driver exchanges at intermediate stops, which requires spatial and temporal synchronization of driver routes. 
We presented a mathematical formulation defined on a time-expanded multi-digraph. Furthermore, we proposed a solution approach that extends an \gls{MIP} implementation by heuristic components and a destructive bound improvement procedure.
We used real-world data provided by one of Europe's leading coach companies to validate our \glsfirst{DBMH}. 
Our solution approach solves instances with up to 390 locations and 70 ride segments to optimality with an average computational time under 210 seconds. Even for large-scale instances we can provide feasible solutions within seconds. In contrast, standalone \gls{MIP} solvers fail to find feasible solutions for a third of the instances even after one hour. 
Compared to standalone \gls{MIP} solvers, our matheuristic improves computational times by 79.42\% for small instances. 
Moreover, compared to current solutions deployed in practice, we can reduce the number of drivers required by 13.07\% on average. Routing and scheduling drivers at the subcontractor level allows utilizing synergies between lines, leading to more efficient solutions. We further demonstrated that the consideration of driver exchanges is beneficial in long-distance bus networks and leads to average savings of 42.69\%. 
Considering driver exchanges at intermediate stops in addition to driver exchanges at regular stops allows for further cost savings, in particular for small and medium instances. For large instances, such a comparison cannot be drawn as many instances cannot be solved to optimality.

Our results show that the instances differ significantly in terms of synchronization potential, which raises a variety of interesting research questions. For example, using machine learning approaches to identify attributes that indicate which instances have a high chance of benefiting from the consideration of driver exchanges at intermediate stops or planning on subcontractor level might be a worthwhile research direction. 
Next, exploiting decomposition-based matheuristics that leverage column generation for subproblems also arising in railway or airline crew scheduling, i.e., problems in which trip legs and departure times are fixed, remains an interesting avenue for future research. Furthermore, adapting existing column generation approaches to consider intermediate stops seems promising. 
Another path for future research is to integrate the \gls{DRSPMS} with other so far separately solved planning problems such as the assignment of rides to buses.
From a methodological and a practical point of view, including uncertainty of travel times and developing a robust solution approach is another important step. 
Finally, studying the potential of driver exchanges en route in other real-world applications such as truck platooning, particularly with increasing automation levels, seems promising.

\section*{Acknowledgments}
We thank Dr. Berit Johannes and colleagues from Flix SE for providing and assisting with the real-world data set, and for providing valuable input and feedback.
The authors gratefully acknowledge the computational and data resources provided by the Leibniz Supercomputing Centre (www.lrz.de).

\section*{Funding}
This work is funded by the Deutsche Forschungsgemeinschaft (DFG, German Research Foundation) – Project Number 277991500.

\section*{Declaration of Interest}
None.


\newpage
\bibliographystyle{ormsv080}
\bibliography{./references}

\begin{thebibliography}{32}
\expandafter\ifx\csname natexlab\endcsname\relax\def\natexlab#1{#1}\fi
\expandafter\ifx\csname url\endcsname\relax
  \def\url#1{{\tt #1}}\fi
\expandafter\ifx\csname urlprefix\endcsname\relax\def\urlprefix{URL }\fi
\expandafter\ifx\csname urlstyle\endcsname\relax
  \expandafter\ifx\csname doi\endcsname\relax
  \def\doi#1{doi:\discretionary{}{}{}#1}\fi \else
  \expandafter\ifx\csname doi\endcsname\relax
  \def\doi{doi:\discretionary{}{}{}\begingroup \urlstyle{rm}\Url}\fi \fi

\bibitem[{Bach et~al.(2016)Bach, Dollevoet, and
  Huisman}]{BachDollevoetEtAl2016}
Bach, Lukas, Twan Dollevoet, Dennis Huisman. 2016.
\newblock Integrating {{Timetabling}} and {{Crew Scheduling}} at a {{Freight
  Railway Operator}}.
\newblock {\it Transportation Science\/} {\bf 50}(3) 878--891.

\bibitem[{Barnhart et~al.(2003)Barnhart, Cohn, Johnson, Klabjan, Nemhauser, and
  Vance}]{BarnhartCohnEtAl2003}
Barnhart, C., A.~M. Cohn, E.~L. Johnson, D.~Klabjan, G.~L. Nemhauser, P.~H.
  Vance. 2003.
\newblock Airline {{Crew Scheduling}}.
\newblock R.~W. Hall, ed., {\it Handbook of {{Transportation Science}}\/}.
  International {{Series}} in {{Operations Research}} \& {{Management
  Science}}, {Springer US}, {Boston, MA}, 517--560.

\bibitem[{Boland et~al.(2017)Boland, Hewitt, Marshall, and
  Savelsbergh}]{BolandHewittEtAl2017}
Boland, N., M.~Hewitt, L.~Marshall, M.~Savelsbergh. 2017.
\newblock The {{Continuous-Time Service Network Design Problem}}.
\newblock {\it Operations Research\/} {\bf 65}(5) 1303--1321.

\bibitem[{Boland et~al.(2019)Boland, Hewitt, Marshall, and
  Savelsbergh}]{BolandHewittEtAl2019}
Boland, N., M.~Hewitt, L.~Marshall, M.~Savelsbergh. 2019.
\newblock The price of discretizing time: A study in service network design.
\newblock {\it EURO Journal on Transportation and Logistics\/} {\bf 8}(2)
  195--216.

\bibitem[{Boland and Savelsbergh(2019)}]{BolandSavelsbergh2019}
Boland, Natashia~L., Martin W.~P. Savelsbergh. 2019.
\newblock Perspectives on integer programming for time-dependent models.
\newblock {\it TOP\/} {\bf 27}(2) 147--173.

\bibitem[{Deveci and Demirel(2018)}]{DeveciDemirel2018}
Deveci, M., N.~{\c C}. Demirel. 2018.
\newblock A survey of the literature on airline crew scheduling.
\newblock {\it Engineering Applications of Artificial Intelligence\/} {\bf 74}
  54--69.

\bibitem[{{Dom{\'i}nguez-Mart{\'i}n} et~al.(2018){Dom{\'i}nguez-Mart{\'i}n},
  {Rodr{\'i}guez-Mart{\'i}n}, and
  {Salazar-Gonz{\'a}lez}}]{Dominguez-MartinRodriguez-MartinEtAl2018}
{Dom{\'i}nguez-Mart{\'i}n}, B., I.~{Rodr{\'i}guez-Mart{\'i}n}, J.-J.
  {Salazar-Gonz{\'a}lez}. 2018.
\newblock The driver and vehicle routing problem.
\newblock {\it Computers \& Operations Research\/} {\bf 92} 56--64.

\bibitem[{Drexl(2012)}]{Drexl2012}
Drexl, M. 2012.
\newblock Synchronization in {{Vehicle Routing}}\textemdash{{A Survey}} of
  {{VRPs}} with {{Multiple Synchronization Constraints}}.
\newblock {\it Transportation Science\/} {\bf 46}(3) 297--316.

\bibitem[{Drexl et~al.(2013)Drexl, Rieck, Sigl, and Press}]{DrexlRieckEtAl2013}
Drexl, M., J.~Rieck, T.~Sigl, B.~Press. 2013.
\newblock Simultaneous {{Vehicle}} and {{Crew Routing}} and {{Scheduling}} for
  {{Partial-}} and {{Full-Load Long-Distance Road Transport}}.
\newblock {\it Business Research\/} {\bf 6}(2) 242--264.

\bibitem[{{European Union}(2006)}]{EuropeanUnion2006}
{European Union}. 2006.
\newblock Regulation ({{EC}}) {{No}} 561/2006 of the {{European Parliament}}
  and of the {{Council}} of 15 {{March}} 2006 on the harmonisation of certain
  social legislation relating to road transport.
\newblock {\it Official Journal of the European Union\/} {\bf 102} 1--13.

\bibitem[{Fink et~al.(2019)Fink, Desaulniers, Frey, Kiermaier, Kolisch, and
  Soumis}]{FinkDesaulniersEtAl2019}
Fink, M., G.~Desaulniers, M.~Frey, F.~Kiermaier, R.~Kolisch, F.~Soumis. 2019.
\newblock Column generation for vehicle routing problems with multiple
  synchronization constraints.
\newblock {\it European Journal of Operational Research\/} {\bf 272}(2)
  699--711.

\bibitem[{Goel and Irnich(2016)}]{GoelIrnich2016}
Goel, A., S.~Irnich. 2016.
\newblock An {{Exact Method}} for {{Vehicle Routing}} and {{Truck Driver
  Scheduling Problems}}.
\newblock {\it Transportation Science\/} {\bf 51}(2) 737--754.

\bibitem[{Guastaroba et~al.(2016)Guastaroba, Speranza, and
  Vigo}]{GuastarobaSperanzaEtAl2016}
Guastaroba, G., M.~G. Speranza, D.~Vigo. 2016.
\newblock Intermediate {{Facilities}} in {{Freight Transportation Planning}}:
  {{A Survey}}.
\newblock {\it Transportation Science\/} {\bf 50}(3) 763--789.

\bibitem[{Heil et~al.(2020)Heil, Hoffmann, and Buscher}]{HeilHoffmannEtAl2020}
Heil, J., K.~Hoffmann, U.~Buscher. 2020.
\newblock Railway crew scheduling: {{Models}}, methods and applications.
\newblock {\it European Journal of Operational Research\/} {\bf 283}(2)
  405--425.

\bibitem[{Hollis et~al.(2006)Hollis, Forbes, and
  Douglas}]{HollisForbesEtAl2006}
Hollis, B.~L., M.~A. Forbes, B.~E. Douglas. 2006.
\newblock Vehicle routing and crew scheduling for metropolitan mail
  distribution at {{Australia Post}}.
\newblock {\it European Journal of Operational Research\/} {\bf 173}(1)
  133--150.

\bibitem[{{Ibarra-Rojas} et~al.(2015){Ibarra-Rojas}, Delgado, Giesen, and
  Mu{\~n}oz}]{Ibarra-RojasDelgadoEtAl2015}
{Ibarra-Rojas}, O.~J., F.~Delgado, R.~Giesen, J.~C. Mu{\~n}oz. 2015.
\newblock Planning, operation, and control of bus transport systems: {{A}}
  literature review.
\newblock {\it Transportation Research Part B: Methodological\/} {\bf 77}
  38--75.

\bibitem[{Kergosien et~al.(2011)Kergosien, Lent{\'e}, Piton, and
  Billaut}]{KergosienLenteEtAl2011}
Kergosien, Y., {\relax Ch}~Lent{\'e}, D.~Piton, J.~C. Billaut. 2011.
\newblock A tabu search heuristic for the dynamic transportation of patients
  between care units.
\newblock {\it European Journal of Operational Research\/} {\bf 214}(2)
  442--452.

\bibitem[{Kim et~al.(2010)Kim, Koo, and Park}]{KimKooEtAl2010}
Kim, B.-I., J.~Koo, J.~Park. 2010.
\newblock The combined manpower-vehicle routing problem for multi-staged
  services.
\newblock {\it Expert Systems with Applications\/} {\bf 37}(12) 8424--8431.

\bibitem[{Klein and Scholl(1999)}]{KleinScholl1999}
Klein, R., A.~Scholl. 1999.
\newblock Computing lower bounds by destructive improvement: {{An}} application
  to resource-constrained project scheduling.
\newblock {\it European Journal of Operational Research\/} {\bf 112}(2)
  322--346.

\bibitem[{Kliewer et~al.(2012)Kliewer, Amberg, and
  Amberg}]{KliewerAmbergEtAl2012}
Kliewer, Natalia, Bastian Amberg, Boris Amberg. 2012.
\newblock Multiple depot vehicle and crew scheduling with time windows for
  scheduled trips.
\newblock {\it Public Transport\/} {\bf 3}(3) 213--244.

\bibitem[{Koub{\^a}a et~al.(2016)Koub{\^a}a, Dhouib, Dhouib, and
  El~Mhamedi}]{KoubaaDhouibEtAl2016}
Koub{\^a}a, M., S.~Dhouib, D.~Dhouib, A.~El~Mhamedi. 2016.
\newblock Truck {{Driver Scheduling Problem}}: {{Literature Review}}.
\newblock {\it IFAC-PapersOnLine\/} {\bf 49}(12) 1950--1955.

\bibitem[{Lam et~al.(2020)Lam, Van~Hentenryck, and
  Kilby}]{LamVanHentenryckEtAl2020}
Lam, E., P.~Van~Hentenryck, P.~Kilby. 2020.
\newblock Joint {{Vehicle}} and {{Crew Routing}} and {{Scheduling}}.
\newblock {\it Transportation Science\/} {\bf 54}(2) 488--511.

\bibitem[{Larrain et~al.(2017)Larrain, Coelho, and
  Cataldo}]{LarrainCoelhoEtAl2017}
Larrain, H., L.~C. Coelho, A.~Cataldo. 2017.
\newblock A {{Variable MIP Neighborhood Descent}} algorithm for managing
  inventory and distribution of cash in automated teller machines.
\newblock {\it Computers \& Operations Research\/} {\bf 85} 22--31.

\bibitem[{Meisel and Kopfer(2014)}]{MeiselKopfer2014}
Meisel, F., H.~Kopfer. 2014.
\newblock Synchronized routing of active and passive means of transport.
\newblock {\it OR Spectrum\/} {\bf 36}(2) 297--322.

\bibitem[{Schiffer et~al.(2017)Schiffer, Laporte, Schneider, and
  Walther}]{SchifferLaporteEtAl2017}
Schiffer, M., G.~Laporte, M.~Schneider, G.~Walther. 2017.
\newblock The impact of synchronizing driver breaks and recharging operations
  for electric vehicles.
\newblock Tech. Rep. G-2017-46, {GERAD, HEC Montreal, Canada}.

\bibitem[{Schiffer et~al.(2019)Schiffer, Schneider, Walther, and
  Laporte}]{SchifferSchneiderEtAl2019}
Schiffer, M., M.~Schneider, G.~Walther, G.~Laporte. 2019.
\newblock Vehicle {{Routing}} and {{Location Routing}} with {{Intermediate
  Stops}}: {{A Review}}.
\newblock {\it Transportation Science\/} {\bf 53}(2) 319--343.

\bibitem[{Schiffer and Walther(2018)}]{SchifferWalther2018}
Schiffer, M., G.~Walther. 2018.
\newblock An {{Adaptive Large Neighborhood Search}} for the {{Location-routing
  Problem}} with {{Intra-route Facilities}}.
\newblock {\it Transportation Science\/} {\bf 52}(2) 331--352.

\bibitem[{Schneider and L{\"o}ffler(2019)}]{SchneiderLoffler2019}
Schneider, M., M.~L{\"o}ffler. 2019.
\newblock Large {{Composite Neighborhoods}} for the {{Capacitated
  Location-Routing Problem}}.
\newblock {\it Transportation Science\/} {\bf 53}(1) 301--318.

\bibitem[{Tilk et~al.(2018)Tilk, Bianchessi, Drexl, Irnich, and
  Meisel}]{TilkBianchessiEtAl2018}
Tilk, C., N.~Bianchessi, M.~Drexl, S.~Irnich, F.~Meisel. 2018.
\newblock Branch-and-{{Price-and-Cut}} for the {{Active-Passive Vehicle-Routing
  Problem}}.
\newblock {\it Transportation Science\/} {\bf 52}(2) 300--319.

\bibitem[{Tilk et~al.(2019)Tilk, Drexl, and Irnich}]{TilkDrexlEtAl2019}
Tilk, C., M.~Drexl, S.~Irnich. 2019.
\newblock Nested branch-and-price-and-cut for vehicle routing problems with
  multiple resource interdependencies.
\newblock {\it European Journal of Operational Research\/} {\bf 276}(2)
  549--565.

\bibitem[{Tilk and Goel(2020)}]{TilkGoel2020}
Tilk, C., A.~Goel. 2020.
\newblock Bidirectional labeling for solving vehicle routing and truck driver
  scheduling problems.
\newblock {\it European Journal of Operational Research\/} {\bf 283}(1)
  108--124.

\bibitem[{Yin et~al.(2021)Yin, D'Ariano, Wang, Yang, and
  Tang}]{YinDArianoEtAl2021}
Yin, Jiateng, Andrea D'Ariano, Yihui Wang, Lixing Yang, Tao Tang. 2021.
\newblock Timetable coordination in a rail transit network with time-dependent
  passenger demand.
\newblock {\it European Journal of Operational Research\/} {\bf 295}(1)
  183--202.

\end{thebibliography}

\newpage

\section*{Appendix A: \gls{MIP} with Consideration of the Subcontractor's Objectives}\label{app:postOpt}
To create driver schedules according to the subcontractor's objective with the \gls{MIP} defined in Section \ref{ss:mip}, we need to modify some constraints and variables. In the following, we describe the modifications needed. 
Specifically, we consider the following five objective functions: i)~minimize the total workload, ii)~balance the workload, iii)~minimize the total steering time, iv)~balance the steering time, v)~minimize the total ride duration.
Given an optimal driver count~$\solX$ resulting from the solution of the original problem as formalized in Section \ref{ss:mip}, we first fix the number of drivers to this value.  We do so by tightening the set of drivers to $\setOfDrivers = \{1, ..., \solX\}$, fixing the lower bound to $LB = \solX$, and replacing Constraint (\ref{eq:LB_x}) with Equality (\ref{eq:LBfixed}):
\begin{align}
\sum_{k \in \setOfDrivers} \sum_{i \in \mathscr{T}(0)}\sum_{(i,j,m) \in \delta^{+}(i)}\varX_{kijm} \; = \; LB \label{eq:LBfixed}
\end{align}

Next, we can integrate the objective function. To minimize the total workload, we replace objective function (\ref{eq:objective}) with (\ref{eq:minTotalWL}):
\begin{align}
\text{min} \qquad \sum_{k \in \setOfDrivers} (\sum_{i \in \mathscr{T}(0)} \sum_{(j,i,m) \in \delta^{-}(i)} \paramNodeTime_{j} \varX_{kjim} - \sum_{i \in \mathscr{T}(0)} \sum_{(i,j,m) \in \delta^{+}(i)}   \paramNodeTime_{j}\varX_{kijm}) \label{eq:minTotalWL}
\end{align}

\noindent
To minimize the total steering time, we replace objective function (\ref{eq:objective}) with (\ref{eq:minTotalSteering}):
\begin{align}
\text{min} \qquad \sum_{k \in \setOfDrivers} \sum_{\substack{(i,j,m) \in \setOfArcsTimeExp:\\ m=1}} \paramConsumption_{ijm} \varX_{kijm} \label{eq:minTotalSteering}
\end{align}

\noindent
To minimize the total duration of all bus rides, we replace objective function (\ref{eq:objective}) with (\ref{eq:minTotalRideDur}):
\begin{align}
\text{min} \qquad \sum_{r \in \setOfRides} \sum_{k \in \setOfDrivers} (
\sum_{\substack{(j,i,m) \in \setOfArcsTimeExp:\\ i \in \mathscr{T}(d_{r}), m=1}}  \paramNodeTime_{i} \varX_{kjim} -
\sum_{\substack{(i,j,m) \in \setOfArcsTimeExp:\\ i \in \mathscr{T}(p_{r}), m=1}}  \paramNodeTime_{i} \varX_{kijm}) 
\label{eq:minTotalRideDur}
\end{align}

\noindent
To integrate balancing objective functions, we define two additional decision variables $\bar{\varB} \geq 0$ and $\underline{\varB} \leq 0$. We also add constraints connecting the balancing variables with the resources to be balanced. For balancing the workload, we replace objective function (\ref{eq:objective}) with (\ref{eq:balancingOF}) and add constraints (\ref{eq:balanceWorkloadUpper}) - (\ref{eq:balanceWorkloadLower}) and variable domains~(\ref{eq:domainBupper}) - (\ref{eq:domainBlower}) to the model:
\begin{align}
\text{min} \qquad \bar{\varB} + \underline{\varB} \label{eq:balancingOF}
\end{align}

\vspace{-0.8cm}
\begin{align}
\sum_{i \in \mathscr{T}(0)} \sum_{(j,i,m) \in \delta^{-}(i)} \paramNodeTime_{j} \varX_{kjim} - \sum_{i \in \mathscr{T}(0)} \sum_{(i,j,m) \in \delta^{+}(i)}   \paramNodeTime_{j}\varX_{kijm} \leq \bar{\varB} \qquad & \forall \; k \in \setOfDrivers \label{eq:balanceWorkloadUpper}\\
-(\sum_{i \in \mathscr{T}(0)} \sum_{(j,i,m) \in \delta^{-}(i)} \paramNodeTime_{j} \varX_{kjim} - \sum_{i \in \mathscr{T}(0)} \sum_{(i,j,m) \in \delta^{+}(i)}   \paramNodeTime_{j}\varX_{kijm}) \leq \underline{\varB} \qquad & \forall \; k \in \setOfDrivers \label{eq:balanceWorkloadLower}\\
\bar{\varB} \geq 0 \qquad &   \label{eq:domainBupper}\\
\underline{\varB} \leq 0 \qquad &   \label{eq:domainBlower}
\end{align}

\noindent
To balance the steering time, we use the same objective function (\ref{eq:balancingOF}) and variable domains (\ref{eq:domainBupper}) - (\ref{eq:domainBlower}), but instead of constraints (\ref{eq:balanceWorkloadUpper}) - (\ref{eq:balanceWorkloadLower}) we add constraints (\ref{eq:balanceSteeringUpper}) - (\ref{eq:balanceSteeringLower}) to the model:
\begin{align}
\sum_{\substack{(i,j,m) \in \setOfArcsTimeExp:\\ m=1}} \paramConsumption_{ijm} \varX_{kijm} \; \leq \; \bar{\varB} \qquad & \forall \; k \in \setOfDrivers \label{eq:balanceSteeringUpper}\\
-\sum_{\substack{(i,j,m) \in \setOfArcsTimeExp:\\ m=1}} \paramConsumption_{ijm} \varX_{kijm} \; \leq \; \underline{\varB} \qquad & \forall \; k \in \setOfDrivers \label{eq:balanceSteeringLower}
\end{align}

Preliminary experiments showed that the problem complexity and improvement potential vary notably depending on the objective function used. 
For a 10\% subset of all small instances, Table \ref{tab:postOpt} depicts the average computational time to optimality (t[s]), the absolute improvement in minutes compared to the optimal solution of the original problem as defined in Sections \ref{ss:mip} and \ref{ss:bounds} ($\Delta \objectiveValue$[min]), and the relative improvement ($\Delta \objectiveValue$[\%]), respectively.
As expected, the improvement potential of the total steering time and the total ride duration is below 1\% and thus negligible. This is due to the detour limit~$\paramDetourLimit$, which allows only minor deviations from the minimum driving times required for a ride. 
On the other hand, the total workload can be reduced by 17.58\% or 247.38 minutes on average. With more than 60\%, the improvement potential for balancing objectives is even higher. However, this optimization potential comes at the cost of an increase in computational time by more than a factor of 100 compared to the \gls{DRSPMS} minimizing the total driver count.

In general, the \gls{MIP} can solve only small instances and is not tractable for medium or large instances. We leave the development of advanced solution methods for the subcontractors' objectives for future research.
\begin{figure}[]
	\begin{floatrow}
		\ttabbox[.7\textwidth]{%
			\small
			\begin{tabular}{lllrrr} 					
				\toprule
				\multicolumn{2}{l}{Objective}&& t[s] & $\Delta \objectiveValue$[min] & $\Delta \objectiveValue$[\%]\\
				\midrule
				i) & Minimize workload && 43.35  & 247.38 & 17.58 \\
				ii) & Balance workload && 24.55 & 153.60 & 65.24\\
				\midrule
				iii) & Minimize steering && 0.12  & 1.42 & 0.17\\
				iv) & Balance steering && 77.24 & 75.99 & 62.01\\
				\midrule
				v) & Minimize ride duration && 0.18 & 2.50 & 0.26\\
				\bottomrule	
			\end{tabular}
		}{%
			\caption{Post-optimization on subcontractor-level}%
			\label{tab:postOpt}%
		}
	\end{floatrow}
\end{figure}

\section*{Appendix B: Computational Analysis of the Lower Bounds}\label{app:lowerBounds}
We conducted a computational analysis based on a subset of instances to evaluate the quality of our different constructive lower bounds described in Section~\ref{ss:bounds}. We computed lower bounds  $LB_1$ and $LB_2$, the combination of both ($LB$), and the linear programming relaxation ($LP$) for every instance in our subset. Table \ref{tab:lbComp} shows the share of instances for which one lower bound dominates the other in terms of objective value. Between $LB_1$ and $LB_2$, no approach can be proven to be superior over the other. 
\begin{figure}[b]
	\begin{floatrow}
		\ttabbox{%
			\small
			\begin{tabular}{lrrrr} 					
				\toprule
				dominates & $LB_1$ & $LB_2$ & $LB$ & $LP$\\
				\midrule
				$LB_1$  & - & 0.6954 & 0.000 & 0.4894 \\
				$LB_2$ & 0.0172 & - & 0.000 & 0.000  \\
				\textbf{\textit{LB}} & \textbf{0.0172} & \textbf{0.6954} & - & \textbf{0.4894} \\
				$LP$ & 0.0071 & 0.0426 & 0.000 & - \\
				\bottomrule	
			\end{tabular}				
		}{%
			\caption{Comparison of lower bounds}%
			\label{tab:lbComp}%
		}
	\end{floatrow}
\vspace{.3cm}
\end{figure}

While bound~$LB_1$ is superior over bound $LB_2$ for 69.54\% of the instances, there also exist instances (1.72\%) for which $LB_2$ is stronger. Furthermore, bound~$LB_1$ outperforms the \gls{LP} relaxation for 48.94\% of the instances. 
The main weakness of the \gls{LP} relaxation is that by allowing fractional arc flows, steering times and working times are not increased by the respective arc's full steering or working time but only by a fraction of it. Thus, none of the working and steering time regulations are sufficiently considered in the \gls{LP} relaxation. Instead, bound~$LB_1$ relaxes mainly the continuous steering time limit. However, bound~$LB_1$ does not consider that a driver can serve only one ride at a time. By combining bound~$LB_1$ with $LB_2$, which accounts for the latter, the combined bound~$LB$ dominates the \gls{LP} relaxation for all instances.
Furthermore, the combined bound~$LB$ dominates bounds~$LB_1$ and $LB_2$ by definition.
Therefore, we use the combined lower bound~$LB$ in our algorithmic framework.

\section*{Appendix C: Computational Analysis of the Algorithmic Components}\label{app:algComponents}
To evaluate the impact of the different components in our solution approach, we run a computational study including 10\% of all instances. Specifically, we consider the following five components: the construction heuristic (\gls{CH}), the \gls{LS} to improve the constructive start solution (\gls{LS}), the destructive bound improvement procedure (\gls{DBI}), callback routines applying the \gls{LS} to new incumbents (\gls{CB}), and the \gls{MIP} (\gls{MIP}).
We tested seven variants, with Variant~0 being the base variant containing all solution components. For Variants~1~to~6, we switched off one component as shown in Table \ref{tab:solComponents}. Note that both Variant~5 and Variant~6 exclude the MIP and the \gls{CB}. For Variant~5, we increase the remaining components' time limit to the global time limit. In contrast, for Variant~6, we keep the time limits of the remaining components unchanged.
For every variant, we conducted 10 runs per instance. Table \ref{tab:solComponents} shows the share of instances solved to feasibility~(solved[\%]), the share of instances solved to optimality (opt.solved[\%]), and the average deviation over 10 runs from the best objective ($\Delta \objectiveValue$[\%]) and from the best optimality gap ($\Delta \; gap$[\%]).  
\begin{figure}[]
	\begin{floatrow}
		\ttabbox{%
			\small
			\begin{tabular}{lrrrrrrrr} 					
				\toprule
				& & \multicolumn{7}{c}{Variants}\\
				\cmidrule{3-9}
				Components & & 0 & 1 & 2 & 3	& 4 & 5 & 6\\
				\midrule
				\gls{CH}		&& \checkmark & \checkmark & & \checkmark & \checkmark & \checkmark & \checkmark\\
				\gls{LS}	&& \checkmark & \checkmark & &  & \checkmark & \checkmark & \checkmark\\
				\gls{DBI} 	&& \checkmark & & \checkmark & \checkmark & \checkmark & \checkmark & \checkmark\\
				\gls{CB}		&& \checkmark & \checkmark & \checkmark & \checkmark &  &  & \\
				\gls{MIP}				&& \checkmark & \checkmark & \checkmark & \checkmark & \checkmark &  & \\		
				\midrule
				solved [\%]		&& \textbf{100.00} & 100.00 & 71.43 & 100.00 & 100.00 & 100.00 & 100.00\\
				opt.solved [\%] && \textbf{60.46} & 60.17 & 60.17 & 59.89 & 60.17 & 60.17 & 57.31\\
				$\Delta \; \objectiveValue$ [\%]	&& \textbf{0.05}& 0.20 & - & 0.24 & 1.21 & 2.16 & 2.16\\
				$\Delta \; gap$ [\%]		&& \textbf{0.51} & 3.18 & - & 0.82 & 1.74 & 1.98 &	3.86\\
				\bottomrule	
			\end{tabular}				
		}{%
			\caption{Contribution of the algorithmic components}%
			\label{tab:solComponents}%
		}
	\end{floatrow}
\end{figure}

The results in Table \ref{tab:solComponents} show that Variant~0 performs best, implying that every component of our solution approach contributes to the solution quality. 
Removing the \gls{DBI}, as done in Variant~1, increases the average deviation from the best optimality gap by more than factor six and further increases the objective function value.
Due to the lack of a heuristic start solution, Variant~2 fails to find any feasible solution for 10 of 35 instances after $\paramTimeLimit$ seconds. Thus, $\Delta \objectiveValue$ and $\Delta gap$ cannot be determined. 
By including the \gls{CH}, feasible solutions can be generated for all instances. However, without the local search component as in Variant~3, the obtained objective values and the obtained optimality gaps for instances not solved optimally suffer.   
Removing \gls{CB}, as in Variant~4, reduces the quality of obtained solutions, particularly for instances not solved optimally. 
Interestingly, removing the \gls{MIP} component does not impact the share of instances solved to optimality any further if the remaining components' time limit is increased as in Variant 5. Without extending the time limit, however, the share of instances solved to optimality decreases by 4.75\%, see Variant 6. In any case, the objective function value and the gap increase.

\section*{Appendix D: Computational Analysis of Neighborhood Search Variants}\label{app:cndVsVND}
We implemented and compared two different search structures: a standard \gls{VND}, applying all operators defined in Section~\ref{sss:vndNeighborhoods} sequentially, and a \glsfirst{LS} (\gls{LS}), exploring the same operators in a composite setting as outlined in Section~\ref{sss:vndNeighborhoods}. 
To compare both neighborhood structures, we run a computational study on a subset including 10\% of all instances presented in Section~\ref{ss:instances}.
Table \ref{tab:cndVsVnd} shows the average deviation of the objective value ($\Delta \objectiveValue$[\%]), the gap ($\Delta gap$[\%]), and the computational time  ($\Delta t$[\%]) from best-known values grouped by instance size.
The results in Table \ref{tab:cndVsVnd} indicate only minor performance differences between the two neighborhood structures. The composite-neighborhood-based \gls{LS} meets all best-known values in terms of objective value and gap. In contrast, results obtained by the \gls{VND} deviate by 0.36\% and 0.48\% from the best-known values for large instances. Due to the more aggressive search behavior, the \gls{LS} outperforms the \gls{VND} in terms of computational time for small- and medium-sized instances. However, the \gls{VND} structure results in shorter computational times for large instances. 
As we restrict computational times for all runs, we are mainly interested in the approach that results in the best objective value. Thus, we chose the composite setting in the local search phase of our matheuristic.

\begin{figure}[b!]
	\begin{floatrow}
		\ttabbox{%
			\small
			\begin{tabular}{llrrrrrrr} 					
				\toprule							
				& & \multicolumn{3}{c}{\gls{VND}} & & \multicolumn{3}{c}{\gls{LS}} \\
				\cmidrule{3-5}
				\cmidrule{7-9}
				Instance size & &  $\Delta \objectiveValue$[\%] & $\Delta gap$[\%] & $\Delta t$[\%] & & $\Delta \objectiveValue$[\%] & $\Delta gap$[\%]&$\Delta t$[\%]\\
				\midrule
				Small && 	0.00 & 0.00 & 50.92 	&& 0.00 & 0.00 & 3.85 \\
				Medium &&  0.00 & 0.00 & 27.33 	&& 0.00 & 0.00 & 16.41 \\
				Large &&  	0.36 & 0.48 & 0.01 	&& 0.00 & 0.00 & 0.15 \\
				\midrule
				All && 	0.09 & 0.12 & 27.70  & & 0.00 & 0.00 & 8.72 \\
				\bottomrule	
			\end{tabular}				
		}{%
			\caption{Performance of different neighborhood structures}%
			\label{tab:cndVsVnd}%
		}
	\end{floatrow}
\end{figure}

\end{document}